\documentclass[12pt,leqno]{amsart}

\makeatletter
\def\thm@space@setup{%
  \thm@preskip=6pt plus 2pt minus 1pt
  \thm@postskip=0pt
}
\makeatother

\usepackage[margin=1in]{geometry}
\usepackage[T1]{fontenc}
\usepackage{lmodern}
\usepackage{microtype}

\usepackage{parskip}

\usepackage{xcolor}

\usepackage{amsmath,amssymb,mathtools}
\makeatletter
\renewenvironment{proof}[1][\proofname]{%
	\par
	\pushQED{\qed}%
	\normalfont
	\topsep6\p@\@plus6\p@\relax
	\trivlist
	\item[\hskip\labelsep\itshape #1\@addpunct{.}]%
	\ignorespaces
}{%
	\par\nobreak\vspace{2pt}% small space before the box
	\popQED%
	\medskip% space after the box
	\endtrivlist
	\@endpefalse
}
\makeatother

\allowdisplaybreaks

\usepackage{etoolbox}

\usepackage{enumitem}
\setlist[enumerate]{label=(\arabic*),leftmargin=*,itemsep=2pt}
\setlist[itemize]{leftmargin=*,itemsep=2pt}

\usepackage{array,booktabs,makecell,multirow}
\usepackage{tabularx}
\usepackage{adjustbox}

\usepackage{graphicx}
\usepackage{float}
\usepackage{comment}
\usepackage{tikz}

\usepackage[unicode]{hyperref}
\hypersetup{
	colorlinks=true,
	linkcolor=blue!50!black,
	citecolor=blue!50!black,
	urlcolor=blue!50!black
}
\usepackage[nameinlink,noabbrev]{cleveref}

\usepackage[numbers,sort&compress]{natbib}

\DeclareMathOperator{\psl}{PSL}
\DeclareMathOperator{\pgl}{PGL}
\DeclareMathOperator{\slg}{SL}
\DeclareMathOperator{\gl}{GL}
\DeclareMathOperator{\tr}{tr}
\DeclareMathOperator{\st}{st}
\DeclareMathOperator{\ord}{ord}

\numberwithin{equation}{section}
\numberwithin{figure}{section}
\numberwithin{table}{section}

\newtheorem{theorem}[equation]{Theorem}
\newtheorem{thm}[equation]{Theorem}      % alias, if you use both
\newtheorem{prop}[equation]{Proposition}
\newtheorem{lem}[equation]{Lemma}
\newtheorem{cor}[equation]{Corollary}

\theoremstyle{definition}

\newtheorem{exm}[equation]{Example}

\theoremstyle{remark}
\newtheorem{rmk}[equation]{Remark}

\theoremstyle{plain}
\newtheorem{mainthm}{Theorem}

\makeatletter
\let\c@figure\c@equation
\let\c@table\c@equation
\makeatother

\usepackage{multirow}

\title[Lattices of Minimal Covolume in $\slg_2$ over Dyadic Local Fields]
{Quaternion Fixed Balls and Lattices of Minimal Covolume in $\slg_2$ over Dyadic Local Fields}

\author{Yukun Du}
\address[Y. Du]{Department of Mathematics, University of Georgia, Athens, GA 30603}
\email{yukun.du@uga.edu}

\date{\today}


\begin{document}

\begin{abstract}
	In this paper, we determine the lattices in $\mathrm{SL}_2(K)$ with minimal covolume for each finite extension $K/\mathbb{Q}_2$. Using the action of $\mathrm{SL}_2(K)$ on its Bruhat--Tits tree, we investigate the shape of the fixed subtrees of finite subgroups, particularly cyclic and dicyclic subgroups of $2$-power order. We determine which larger vertex stabilizers can occur at boundary vertices of the fixed subtrees of dicyclic subgroups and subgroups isomorphic to $2.A_4$. We use these results to construct graphs of groups defining lattices in $\mathrm{SL}_2(K)$ and prove that these lattices attain the minimal covolume. This extends a result of Lubotzky and Weigel, which treats finite extensions of $\mathbb{Q}_p$ for odd primes $p$.
\end{abstract}
\maketitle
% \tableofcontents
		
\section{Introduction}\label{sec:introduction}
It is a classical problem to determine how \emph{dense} a lattice $\Gamma$ in a Lie group $G$ can be. More precisely, one seeks the minimum of the normalized Haar measure $\mu(\Gamma\backslash G)$. When $G = \slg(2,\mathbb{R})$, a result of Siegel \cite{siegel1945some} implies that the minimum is attained when $\Gamma$ is the preimage of the $(2,3,7)$-triangle group, where $\min_{\Gamma}\mu(\Gamma\backslash G) = \pi/21$. The complex case $G = \slg(2,\mathbb{C})$ was addressed more recently by Gehring--Martin \cite{gehring2009minimal} and Marshall--Martin \cite{marshall2012minimal}. Namely, they determined that
\[
\min_{\Gamma}\mu(\Gamma\backslash \slg(2,\mathbb{C})) = \frac{275^{3/2}}{128\pi^6}\zeta_k(2)\approx 0.03905\dots,
\]
where $k = \mathbb{Q}(\alpha)$ and $\alpha$ is a root of $\alpha^4 +6\alpha^3+12\alpha^2+9\alpha+1 = 0$. This minimal volume is attained by an $\mathbf{H}^3$-orbifold, admitting $60$-sheeted manifold covers, given by either the $-4/3$-Dehn surgery along the $(-2,3,7)$ pretzel knot or the $(-7/2,-10/3)$-Dehn surgery along the Whitehead link \cite{conder2006maximal}. By contrast, the minimal covolume in certain families of higher-rank Lie groups is attained by the standard integral lattices. In particular, Thilmany \cite{thilmany2019lattices} proved that the minimal covolume lattice in $\slg(n,\mathbb{R})$ is $\slg(n,\mathbb{Z})$ for $n\geq 3$; for $\operatorname{Sp}(2n,\mathbb{R})$ with $n\geq 2$, D\v{z}ambi{\'c}--Holm--K\"ohl \cite{dvzambic2024siegel} showed that the lattice of minimal covolume is $\operatorname{Sp}(2n,\mathbb{Z})$.

An analogous question arises for Lie groups over non-archimedean local fields. Lubotzky \cite{lubotzky1990lattices} studied the minimal covolume lattices in $\slg(2,K)$, where $K$ is either a field of positive characteristic $\mathbb{F}_q((t))$ or the $p$-adic field $\mathbb{Q}_p$. Notably, the minimal covolume lattices in the positive-characteristic cases are non-uniform lattices $\slg_2(\mathbb{F}_q[t^{-1}])$. In the higher-rank case, Allcock--Kato \cite{allcock2013densest} determined the lattices of minimal covolume in $\pgl(3,\mathbb{Q}_2)$, with covolume $1/21$. The two minimizing lattices are commensurable, arising from a Hermitian form over $\mathbb{Q}(\sqrt{-7})$.

The case where $K$ is a finite extension of $\mathbb{Q}_p$ is more interesting. When $p$ is an odd prime, Lubotzky--Weigel \cite{lubotzky1999lattices} determined the lattices of minimal covolume in $\slg(2,K)$. The lattices they constructed are described by graphs of finite groups with up to six vertices, whereas paths of length $1$ and loops of length $2$ are enough when $K = \mathbb{Q}_p$.

A few questions arise when generalizing Lubotzky and Weigel's method to the $p=2$ case. For one thing, when $p=2$, elliptic elements in $\slg(2,K)$ of order $p^k$ have more intricate dynamics; in particular, the diameter of the fixed subtree has a more complicated description. Furthermore, many vertex groups of lattices in $\slg(2,K)$ are dicyclic, and the elements of order $4$ in these groups result in more complicated behavior when $p=2$.

We aim to resolve these issues and settle the problem for $\slg_2$ over finite extensions of $\mathbb{Q}_2$. The main result we have obtained is the following:
\begin{mainthm}\label{thm:main}
    Let $K$ be a finite extension of $\mathbb{Q}_2$, of degree $\lvert K/\mathbb{Q}_2\rvert = n$, residue field order $q=2^f$, and ramification index $e = n/f$. Let $\xi_m$ be a primitive $m$-th root of unity, and put $z_m = \xi_m+\xi_m^{-1}$. The Haar measure $\mu$ of $G = \slg(2,K)$ is normalized so that $\mu(\slg(2,\mathcal{O}_K)) = 1$. Then the minimal covolume among cocompact lattices $\Gamma<\slg(2,K)$ is given by the following:
    \[
    \lambda(K) = \begin{cases}
        \frac{1}{120}, & \ K = \mathbb{Q}_2(\xi_3)\text{ (cf. Example \ref{exm:3})},\text{ or } K = \mathbb{Q}_2(\xi_{15})\text{ (cf. Example \ref{exm:4})}, \\
        \frac{1}{48}, & \ K = \mathbb{Q}_2(z_8)\text{ (cf. Example \ref{exm:a})}, \\
        \frac{7}{240}, & \ K = \mathbb{Q}_2(\xi_3,z_8)\text{ (cf. Example \ref{exm:b})}, \\
        \frac{1}{24}, & \ K = \mathbb{Q}_2(z_{12})\text{ (cf. Example \ref{exm:5})},\text{ or }K = \mathbb{Q}_2(\xi_3,z_{20})\text{ (cf. Example \ref{exm:7})},\\
        \frac{5}{96}, & \ K = \mathbb{Q}_2(z_{16})\text{ (cf. Example \ref{exm:d})},\\
        \frac{1}{16}, & \ K = \mathbb{Q}_2(z_8,z_{12})\text{ (cf. Example \ref{exm:c})}, \\
        \frac{1}{8}, & \ K=\mathbb{Q}_2(\sqrt{6}),\mathbb{Q}_2(\sqrt{10}),\text{ or }\mathbb{Q}_2(\sqrt{14})\text{ (cf. Example \ref{exm:6})},\\
        \frac{1}{10}, & \ K = \mathbb{Q}_2(\xi_{12})\text{ (cf. Example \ref{exm:reg}, or Example \ref{exm:y})}, \\
        \frac{1}{6}, & \ K = \mathbb{Q}_2(\xi_4)\text{ or } K = \mathbb{Q}_2(\xi_8) \text{ (cf. Example \ref{exm:reg}, or}\\
        & \text{Examples \ref{exm:x} and \ref{exm:z} respectively)}, \\
        \frac{1}{4(q+1)}, & \ K = \mathbb{Q}_2(\xi_{q-1}),\ q\ne 4,\ 16\text{ (cf. Example \ref{exm:2})}, \\
        \frac{1}{2(q+1)}, & \ \text{in all other cases}\text{ (cf. Example \ref{exm:reg}).}
    \end{cases}
    \]
\end{mainthm}
\begin{cor}\label{cor:maincor}
    There are up to two types of graphs of groups for lattices of minimal covolume in $G = \slg(2,K)$. In addition, such graph of groups is unique except for the following cases:
    \begin{itemize}
        \item $K = \mathbb{Q}_2(\xi_4)$ (cf. Example \ref{exm:reg} and Example \ref{exm:x}); in this case one has $\lambda(K) = \frac{1}{6}$.
        \item $K = \mathbb{Q}_2(\xi_8)$ (cf. Example \ref{exm:reg} and Example \ref{exm:z}); in this case one has $\lambda(K) = \frac{1}{6}$.
        \item $K = \mathbb{Q}_2(\xi_{12})$ (cf. Example \ref{exm:reg} and Example \ref{exm:y}); in this case one has $\lambda(K) = \frac{1}{10}$.
    \end{itemize}
\end{cor}
The basics of the question and a Brauer group criterion for the existence of $2.A_4$ subgroups in $\slg(2,K)$ are discussed in Section \ref{sec:prelim}. In Section \ref{sec:fixed_subtree}, we generalize \cite[Theorem 4.2]{lubotzky1999lattices} to give a complete description of the fixed subtrees of elements in $\slg(2,K)$ of order $2^k$. In Section \ref{sec:vertex_stabilizer}, we discuss the possible vertex groups of $\slg(2,K)$ lattices, including the dicyclic and binary polyhedral ones. We begin Section \ref{sec:construction} by recalling a result on the orbits of vertex stabilizers, followed by the construction of the candidate lattices of minimal covolume. Finally, we complete the minimal covolume classification in Section \ref{sec:proof}, treating separately lattices with and without elements of order $8$. These results, namely Proposition \ref{prop:no_8} and Theorem \ref{thm:order_8}, together complete the proof of Theorem \ref{thm:main}.

\textbf{Acknowledgement.} The author would like to thank Saar Hersonsky for introducing them to the minimal covolume problem for groups over non-archimedean local fields, particularly in the setting of group actions on trees. This introduction led them to the questions explored in this paper.

\section{Preliminaries}\label{sec:prelim}
\subsection{Notations and Basic Facts}
In this paper, unless otherwise stated, we let $K$ be a finite extension of $\mathbb{Q}_2$, with uniformizer $\pi = \pi_K$ and ring of integers $\mathcal{O} = \mathcal{O}_K$. The residue field $\mathcal{O}_K/\pi_K\mathcal{O}_K$ is isomorphic to a finite field $\mathbb{F}_q$, where $q = 2^f$, where $f = f(K/\mathbb{Q}_2)$ is called the \emph{residue degree} of $K/\mathbb{Q}_2$. With the \emph{valuation} normalized by $v_K(\pi_K) = 1$, $e = e(K/\mathbb{Q}_2) = v_K(2)$ is called the \emph{ramification degree}. These quantities satisfy
\[
[K:\mathbb{Q}_2] = e(K/\mathbb{Q}_2)\cdot f(K/\mathbb{Q}_2).
\]
The \emph{Bruhat-Tits tree} of the non-archimedean Lie group $G = \slg(2,K)$ is the quotient $T = \pgl(2,K)/\pgl(2,\mathcal{O}_K)$, which is realized as a $(q+1)$-regular tree with edge lengths normalized to $1$. For any vertex $x\in\mathcal{V}(T)$, the \emph{vertex stabilizer} $G_x$ is isomorphic to $\slg(2,\mathcal{O}_K)$. For any element $g\in \slg(2,K)$ or subgroup $A<\slg(2,K)$, we denote its \emph{fixed subtree} by $T_g$, or $T_A$, respectively; we note that $T_g = T_{\langle g\rangle}$.

For the vertex groups appearing in this paper, we denote the dicyclic group of order $4m$ by $C_{2m}.2$, which is a double cover of the dihedral group $D_{2m}$. We denote the three binary polyhedral groups by $2.A_4$, $2.S_4$ and $2.A_5$, which are double covers of the polyhedral groups $A_4$, $S_4$, and $A_5$, respectively.

In the algebraic closure $\overline{\mathbb{Q}}_2$, we let $\xi_m$ be the $m$-th root of unity, and $z_m = \xi_m+\xi_m^{-1}$. In particular, $\xi_3 = \frac{-1+\sqrt{-3}}{2}$, $\xi_4 = \sqrt{-1}$, $z_5 = \frac{-1+\sqrt{5}}{2}$, $z_8 = \sqrt{2}$, and $z_{12} = \sqrt{3}$.

We use the same normalization of Haar measure on $G$, namely $\mu(\slg(2,\mathcal{O}_K)) = 1$, and the covolume of a lattice is calculated as
\[
\mu(G/\Gamma) = \frac{1}{2}\sum_{x\in\mathcal{V}(T//\Gamma)}\frac{1}{\lvert \Gamma_x\rvert}.
\]
Equivalently, since $2\mid d(x,g.x)$ for any $x\in\mathcal{V}(T)$ and $g\in\slg(2,K)$, the $\slg(2,K)$-action preserves the \emph{red vertices} and \emph{blue vertices} of $T$; the covolume $\mu(G/\Gamma)$ equals the sum of the reciprocals of the stabilizer orders over either all blue vertices or all red vertices. 

The finite subgroups of $\slg(2,K)$ fall into two infinite families of cyclic and dicyclic groups, and three exceptional binary polyhedral ones. Namely:
\begin{thm}[{\cite[Theorem 3.1]{lubotzky1999lattices}}]\label{thm:finite_subgroup}
    A finite subgroup $S<\slg(2,K)$ satisfies one of the following:
    \begin{itemize}
    \item $S$ is isomorphic to a subgroup of $\mu_K.2$, where $\mu_K$ is the group of roots of unity in $K$.
    \item $S$ is isomorphic to a subgroup of $\mu^*_{E/K}.2$ for some quadratic extension $E/K$, where
    \[
    \mu^*_{E/K}=\{x\in\mu_E:N_{E/K}(x)=1\}.
    \]
    \item $S$ is isomorphic to one of $2.A_4$, $2.S_4$, or $2.A_5$.
    \end{itemize}
\end{thm}
\begin{rmk}
    For $q=2^{f(K/\mathbb{Q}_2)}$, $k = \max\{a:\xi_{2^a}\in K\}$ and $k' = \max\{b:\xi_{2^b}\in E\}$, we have $\mu_K.2\cong C_{(q-1)2^k}.2$, $\mu^*_{E/K}.2\cong C_{(q+1)2^{k'}}.2$ if $E/K$ is not ramified, and $\mu^*_{E/K}.2\cong C_{2^{k'}}.2$ if $E/K$ is ramified. Moreover, the dicyclic group $C_{2^k}.2$ is also identified as the \emph{generalized quaternion group} $Q_{2^{k+1}}$.
\end{rmk}
For an element $g$ of odd order, its fixed subtree $T_g$ is either a single vertex or a bi-infinite geodesic:
\begin{lem}[{\cite[Lemma 4.1 (c), the case $p=2$]{lubotzky1999lattices}}]\label{lem:odd_fixed}
    Let $K$ be a finite extension of $\mathbb{Q}_2$, $g\in \slg(2,K)$ is an element of odd order $m>1$, then either $m\mid (q-1)$ or $m\mid (q+1)$. In addition:
    \begin{itemize}
        \item If $m\mid (q+1)$, then the fixed tree $T_g$ is a single vertex.
        \item If $m\mid (q-1)$, then the fixed tree $T_g$ is a bi-infinite path.
    \end{itemize}
\end{lem}

For an element $g$ of order $2^k$, its fixed subtree $T_g$ is related to the ramification number of $K$:
\begin{thm}[{\cite[Theorem 4.2]{lubotzky1999lattices}}, the case $p=2$]\label{thm:regular_width}
    Let $K$ be a finite extension of $\mathbb{Q}_2$, and let $g\in \slg(2,K)$ have order $2^k$, with $k\geq 2$. It follows that $g$ has the roots of unity $\xi_{2^k}$ and $\xi_{2^k}^{-1}$ as eigenvalues, so $z_{2^k} \coloneqq \xi_{2^k} + \xi_{2^k}^{-1}\in K$. Define the ramification number $r_k = e(K/\mathbb{Q}_2(z_{2^k}))$ and put $C = C_{\gl(2,K)}(g)$. Then the fixed subtree is $T_g = N(\mathfrak{C}_g,d_k)$, where $\mathfrak{C}_g$ is the characteristic subtree of $C$, called the \emph{core} of $T_g$\footnote{Lubotzky and Weigel called it the \emph{centre}.}, and $d_k\in\mathbb{N}$ is called the \emph{width} of $T_g$. Moreover, one of the following holds.
    \begin{itemize}
        \item If the root of unity $\xi_{2^k}\in K$, then $C$ is a split torus, $\mathfrak{C}_g$ is a bi-infinite path, and $d_k = r_k$.
        \item If $\xi_{2^k}\notin K$, then $C$ is a non-split torus, and $\lvert K(\xi_{2^k})/K\rvert = 2$:
        \begin{itemize}
            \item If $e(K(\xi_{2^k})/K) = 1$, i.e., the extension is unramified, then $\mathfrak{C}_g$ is a single vertex, and $d_k = r_k$.
            \item If $e(K(\xi_{2^k})/K) = 2$, i.e., the extension is totally ramified, then $\mathfrak{C}_g$ is a single edge, and $d_k\geq \frac{r_k-1}{2}$.
        \end{itemize}
    \end{itemize}
\end{thm}
The following conjugacy class uniqueness for finite subgroups is helpful for constructing candidate graph of groups with desired vertex groups incident with cyclic edge groups.
\begin{prop}[{\cite[Proposition 4.9 (2)]{lubotzky1999lattices}}]\label{prop:cyclic_conj}
    Let $A = \langle g\rangle\cong C_m$ and $x\in \partial T_A$. Then all subgroups $B<G_x$ isomorphic to $C_m$ with $x\in \partial T_B$ belong to the same $\gl(2,K)_x$-conjugacy class.
\end{prop}
Below is a particular example of graph of groups representing a lattice in $\slg(2,K)$ which gives a uniform upper bound for the minimal covolume $\lambda(K)$.
\begin{exm}[{\cite[Example (6.1)]{lubotzky1999lattices}}]\label{exm:reg}
    Let $K$ be a finite extension of $\mathbb{Q}_2$, with residue field order $q = 2^{f(K/\mathbb{Q}_2)}$. Then there exists a lattice $\Gamma<G = \slg(2,K)$ with graph of groups
    \begin{center}
        \begin{tikzpicture}[scale=2.5,every node/.style={font=\footnotesize}]
            \draw[thick] (0,0) --node[above]{$C_2$} (1,0);
            \fill[] (0,0) circle(0.025) node[left]{$C_{2(q+1)}$};
            \fill[] (1,0) circle(0.025) node[right]{$C_{2(q+1)}$};
        \end{tikzpicture}
    \end{center}
    The covolume of this lattice is $\frac{1}{2(q+1)}$.
\end{exm}

\subsection{Quaternion Algebras and a Norm Criterion}
Let $K$ be a finite extension of $\mathbb{Q}_p$. By the Artin--Wedderburn theorem, every central simple
$K$-algebra $A$ is isomorphic to $M_r(D)$, where $r\geq 1$ is uniquely determined and the central division $K$-algebra $D$ is unique up to $K$-algebra isomorphism. Two central simple $K$-algebras $A$ and $A'$ are \emph{Brauer equivalent} if the corresponding division algebras satisfy $D\cong D'$ as $K$-algebras. The equivalence classes, equipped with the operation
\[
[A]+[A']=[A\otimes_K A'],
\]
with neutral element $[K]$, form the \emph{Brauer group} of $K$, denoted by $\operatorname{Br}(K)$. Thus $[A]=0$ if and only if $A\cong M_r(K)$ for some $r\geq 1$, i.e., if $A$ is split. There is a canonical isomorphism:
\[
\operatorname{inv}_K:\ \operatorname{Br}(K)\xrightarrow{\sim}\mathbb Q/\mathbb Z,
\]
and $\operatorname{inv}_K([A])$ is called the \emph{Hasse invariant} of the Brauer class $[A]$.

For every finite extension $E/K$, the Brauer restriction map $\operatorname{res}_{E/K}$: $[A]\mapsto [A\otimes_K E]$ satisfies a degree multiplication relation (see \cite[Chapter XIII, \S3]{serre2013local}):
\[
\operatorname{inv}_E([A\otimes_K E])
=\lvert E/K\rvert\cdot\operatorname{inv}_K([A]).
\]
In particular, for a quaternion algebra $A = (a,b)_K$, $A\otimes_K E$ always splits over $E = K(\sqrt{a})$. Since $\lvert E/K\rvert\in\{1,2\}$, one has $\operatorname{inv}_K([A])\in \{[0],[\frac{1}{2}]\}$.

\begin{lem}\label{lem:square_sum}
    Let $K$ be a finite extension of $\mathbb{Q}_2$. Then the following are equivalent:
    \begin{itemize}
        \item There exist $x,y\in K$ such that $x^2 + y^2 = -1$.
        \item There exist $x,y\in \mathcal{O}_K$ such that $x^2+y^2 = -1$.
        \item $\lvert K/\mathbb{Q}_2\rvert$ is even.
    \end{itemize}
\end{lem}
\begin{proof}
    The second condition clearly implies the first. Now, if $x^2 + y^2 = -1$ holds in $K$, with $v(x)\leq v(y)$, $v(x)<0$, then 
    \[
    \left(\frac{y}{x}\right)^2 + \left(\frac{1}{x}\right)^2 = -1,\ \frac{y}{x},\ \frac{1}{x}\in\mathcal{O}_K,
    \]
    so a solution in $\mathcal{O}_K$ exists as well.

    To show that the third condition is equivalent to them, we consider the quaternion algebra $A = (-1,-1)_{\mathbb{Q}_2}$. The splitting criterion for quaternion algebras (see \cite[Chapter III, \S2]{lam2005introduction}) implies that $A\otimes_{\mathbb Q_2}K\simeq M_2(K)$ if and only if there exist $x,y\in K$ such that $x^2+y^2 = -1$. Since there are no such solutions with $x,y\in\mathbb{Q}_2$, our discussion above gives
    \[
    \operatorname{inv}_{\mathbb{Q}_2}([A]) = [\tfrac{1}{2}].
    \]
    For a finite extension $K/\mathbb{Q}_2$, the degree multiplication then gives
    \[
    \operatorname{inv}_K([A\otimes_{\mathbb{Q}_2}K]) = \lvert K/\mathbb{Q}_2\rvert\cdot \operatorname{inv}_{\mathbb{Q}_2}([A]) = [\tfrac{\lvert K/\mathbb{Q}_2\rvert}{2}],
    \]
    which vanishes if and only if $\lvert K/\mathbb{Q}_2\rvert$ is even. Therefore, 
    \[
    \exists x,y\in K,\  x^2 + y^2 = -1\iff 2\mid \lvert K/\mathbb{Q}_2\rvert.
    \]
\end{proof}

\section{Fixed Subtrees of \texorpdfstring{$2$}{2}-Elements}\label{sec:fixed_subtree}
Lubotzky and Weigel did not further analyze the behavior of the width in case (2b), especially when $r_k$ is even. Nevertheless, we obtain the following result.
\begin{thm}\label{thm:width_formula}
    Let $K$, $g$, $r_k$ and $C$ be defined as above, and suppose $\xi_{2^k}\notin K$ and $e(K(\xi_{2^k})/K) = 2$. Then the width $d_k$ of the fixed tree $T_g = N(T_C, d_k)$ satisfies the following.
    \begin{itemize}
        \item If $r_k$ is odd, then $d_k = \frac{r_k-1}{2}$.
        \item If $r_k$ is even, then $\frac{r_k}{2}\leq d_k\leq r_k-1$.
        % , and each admissible $d_k$ is attained by some extension $K$.
    \end{itemize}
    Also, suppose $\kappa$ is the largest integer such that $\slg(2,K)$ contains an element of order $2^{\kappa}$. Then for $3\leq k\leq \kappa$, we have
    \[
    d_{k-1} = d_k+r_k.
    \]
    In particular, if $r_{\kappa}$ is odd, then
    \[
    d_k = r_k - \frac{r_{\kappa}+1}{2}.
    \]
\end{thm}
We also give a complete and explicit description of the possible values of $d_k$ in Appendix \ref{append:a}. This description is not needed for the proof of Theorem \ref{thm:main}.
\begin{lem}\label{lem:max_valuation}
    Put $z=z_{2^k}$ and $w=2-z$. Under the assumptions of Theorem \ref{thm:width_formula}, we have
    \[
    2d_k + 1 = \max_{\alpha\in \mathcal{O}_K}\frac{v(\alpha^2 - z\alpha + 1)}{v(\pi)} = \max_{\beta\in \mathcal{O}_K}\frac{v(\beta^2 + w\beta + w)}{v(\pi)}.
    \]
\end{lem}
\begin{proof}
    The group $\psl(2,K)$ acts on pairs of endpoints (and $\pgl(2,K)$ on triples) in $\partial_\infty T_K = \mathbf{P}^1(K)$. Therefore, up to a suitable $\psl(2,K)$-conjugation of $g$, we may assume that a longest path in $T_g$ lies in the geodesic connecting $0$ and $\infty\in \mathbf{P}^1(K)$. On the one hand, the length of this path equals $2d_k+1$. On the other hand, since $\ord(g) = 2^k$, it follows that $\tr(g) = z_{2^k}$. Therefore we can assume
    \[
    g = \begin{pmatrix}
        \alpha & x \\ x^{-1}(\alpha(z-\alpha)-1) & z-\alpha
    \end{pmatrix},
    \]
    for some $\alpha$ and $x\in K$. Since the stabilizers of the vertices on the path $(0,\infty)$ are
    \[
    \left\{\left.\begin{psmallmatrix}
        a & b \\ c & d
    \end{psmallmatrix}\in\slg(2,K)\right|a,d\in \mathcal{O},\ b\in \pi^i\mathcal{O},\ c\in\pi^{-i}\mathcal{O}\right\},\quad i\in\mathbb{Z},
    \]
    it follows that $\alpha\in \mathcal{O}_K$, and
    \[
    \ell(T_g\cap (0,\infty)) = \frac{v(x) + v(x^{-1}(\alpha(z-\alpha)-1)}{v(\pi)} = \frac{v(\alpha^2-z\alpha+1)}{v(\pi)}.
    \]
    The reverse inequality is immediate. For the other equivalent expression, we note that
    \[
    \alpha^2 - z\alpha + 1 = \alpha^2 - (2-w)\alpha + 1 = (\alpha-1)^2 + w(\alpha - 1) + w,
    \]
    and $\beta = \alpha - 1\in\mathcal{O}_K$ if and only if $\alpha\in\mathcal{O}_K$.
\end{proof}
\begin{lem}[cf. {\cite[Proposition 4.4]{lubotzky1999lattices}}]\label{lem:width_ineq}
    We have
    \[
    r_kv(\pi)\leq \max_{\alpha\in \mathcal{O}_K}v(\alpha^2 - z\alpha + 1)\leq 2r_kv(\pi).
    \]
\end{lem}
\begin{proof}
    For the first inequality, we take $\alpha = 1$:
    \[
    \max_{\alpha\in \mathcal{O}_K}v(\alpha^2 - z\alpha + 1)\geq v(2-z) = v(w) = r_k v(\pi).
    \]
    For the second inequality, we note that in the field extension $K(\xi_{2^k})$,
    \[
    (\alpha^2 - z\alpha + 1) = (\alpha - \xi)(\alpha - \xi^{-1}),
    \]
    so
    \[
    v(\alpha^2 - z\alpha + 1) = v(\alpha - \xi)+v(\alpha - \xi^{-1}).
    \]
    Since $\alpha\in K$, the generator $\sigma\in \operatorname{Gal}(K(\xi_{2^k})/K)$ satisfies $\sigma(\alpha - \xi) = \alpha - \xi^{-1}$, so
    \[
    v(\alpha - \xi)+v(\alpha - \xi^{-1}) = 2v(\alpha - \xi).
    \]
    By the non-archimedean triangle inequality,
    \[
    v(\alpha - \xi) = \min(v(\alpha - \xi),v(\alpha - \xi^{-1}))\leq v((\alpha - \xi^{-1})-(\alpha - \xi)) = v(\xi - \xi^{-1}) = r_kv(\pi).
    \]
    Since the inequality holds for every $\alpha\in\mathcal{O}_K\subset K$, we have
    \[
    \max_{\alpha\in \mathcal{O}_K}v(\alpha^2 - z\alpha + 1)\leq 2r_kv(\pi).
    \]
\end{proof}
\begin{rmk}
    If $K$ is a finite extension of $\mathbb{Q}_p$ for an odd prime $p$, the same arguments give
    \[
    \max_{\alpha\in K}v(\alpha^2 - z\alpha + 1)\leq 2v(\xi_{p^k} - \xi_{p^k}^{-1}) = 2v(\xi_{p^k}-1) = r_kv(\pi),
    \]
    which differs from the case $p=2$ by a factor of $2$ (cf. \cite[proof of Theorem 4.2]{lubotzky1999lattices}). Combining this with the other inequality,
    \[
    \max_{\alpha\in K}v(\alpha^2 - z\alpha + 1)\geq v(2-z_{p^k}) = v((1-\xi_{p^k})^2\xi_{p^k}^{-1}) = r_kv(\pi),
    \]
    we obtain the equality $\max_{\alpha\in K}v(\alpha^2 - z\alpha + 1) = r_kv(\pi)$, leading to $d_k = \frac{r_k-1}{2}$. This agrees with the formula in \cite{lubotzky1999lattices}.
\end{rmk}
\begin{proof}[Proof of Theorem \ref{thm:width_formula}]
    By Lemma \ref{lem:max_valuation}, it suffices to compute $\max_{\beta\in\mathcal{O}_K}v(\beta^2 + w\beta + w)$ for showing the prescribed range of $d_k$. If $r_k$ is odd, then
    \[
    v(\beta^2) = 2v(\beta)\neq r_kv(\pi) = v(w).
    \]
    Therefore,
    \[
    v(\beta^2 + w\beta + w) = \min(v(\beta^2),v(w))\leq r_k v(\pi),
    \]
    and equality holds when $v(\beta)>v(w)/2$, e.g., when $\beta = 0$. Therefore
    \[
    2d_k+1 = \max_{\beta\in\mathcal{O}_K}\frac{v(\beta^2 + w\beta + w)}{v(\pi)} = r_k,
    \]
    and $d_k = \frac{r_k-1}{2}$.

    If $r_k$ is even, by Lemma \ref{lem:width_ineq},
    \[
    r_k\leq 2d_k+1\leq 2r_k,
    \]
    so
    \[
    \frac{r_k-1}{2}\leq d_k\leq \frac{2r_k-1}{2}.
    \]
    Since $d_k$ is an integer, we further have
    \[
    \frac{r_k}{2}\leq d_k\leq r_k-1.
    \]
    For the width difference relationship, recall that the proofs of Lemmas \ref{lem:max_valuation} and \ref{lem:width_ineq} give
    \[
    2d_k+1 = \max_{\alpha\in\mathcal{O}_K} 2v_K(\alpha - \xi_{2^k}) = \max_{\alpha\in\mathcal{O}_K}v_L(\alpha - \xi_{2^k}),
    \]
    where $L = K(i) = K(\xi_{2^{\kappa}})$. Fix a uniformizer $\pi_L$, and write $\xi_{2^k} = a_k + b_k\pi_L$, where $a_k,b_k\in K$. Therefore,
    \[
    v_L(\alpha - \xi_{2^k}) = \min\{v_L(\alpha - a_k),v_L(b_k\pi_L)\} = \min\{2v_K(\alpha - a_k),2v_K(b_k)+1\}.
    \]
    Comparing parities, we see that $d_k = v_K(b_k)$. Applying the nontrivial automorphism $\sigma\in\operatorname{Aut}(L/K)$, we eliminate $a_k$ and obtain
    \[
    b_k = \frac{\xi_{2^k} - \xi_{2^k}^{-1}}{\pi_L - \sigma(\pi_L)}.
    \]
    Combining this with the identity $\xi_{2^{k-1}} - \xi_{2^{k-1}}^{-1} = z_{2^k}(\xi_{2^k} - \xi_{2^k}^{-1})$ gives $b_{k-1} = z_{2^k}b_k$, so
    \[
    d_{k-1} = d_k + v_K(z_{2^k}) = d_k+r_k.
    \]
    The formula when $2\nmid r_{\kappa}$ then follows inductively from $d_{\kappa} = \frac{r_{\kappa}-1}{2}$.
\end{proof}
% As $k$ varies, the widths $d_k$ satisfy a more precise relationship.
% \begin{prop}\label{prop:width_diff_inv}
%     In the setting of Theorem \ref{thm:width_formula}, 
% \end{prop}
% \begin{proof}
    
% \end{proof}
\section{Vertex Stabilizers in \texorpdfstring{$\slg(2,K)$}{PSL(2,K)}-lattices}\label{sec:vertex_stabilizer}
\begin{prop}[cf. {\cite[Proposition 3.3]{lubotzky1999lattices}}]
    Let $K$ be a finite extension of $\mathbb{Q}_2$ with residue field $\mathbb{F}_q$.
    \begin{itemize}
        \item The group $2.A_4$ is isomorphic to a subgroup of $\slg(2, K)$ if and only if $\lvert K/\mathbb{Q}_2\rvert$ is even.
        \item The group $2.S_4$ is isomorphic to a subgroup of $\slg(2, K)$ if and only if $z_8\in K$.
        \item The group $2.A_5$ is isomorphic to a subgroup of $\slg(2, K)$ if and only if $q\equiv \pm 1 (\bmod 5)$, i.e., if $z_5\in K$.
    \end{itemize}
\end{prop}
Note that both conditions (2) and (3) imply (1): both $z_8 = \sqrt{2}$ and $z_5 = \frac{-1+\sqrt{5}}{2}$ are quadratic over $\mathbb{Q}_2$, so extensions containing either element have even degree over $\mathbb{Q}_2$.
\begin{proof}
    The same argument as in the proof of \cite[Proposition 3.3]{lubotzky1999lattices} shows that $2.A_4$ or $2.A_5$ is contained in $\slg(2,K)$ if and only if $A_4$ or $A_5< \pgl(2,K)$.

    (a) We know $A_4<\pgl(2,K)$ if and only if $\exists y,z\in K$ with $y^2 + z^2 = -1$. By Lemma \ref{lem:square_sum}, this is equivalent to the condition that $\lvert K/\mathbb{Q}_2\rvert$ is even.

    (b) We have $2.S_4< \slg(2,K)$ if and only if $S_4 < \psl(2,K)$. On the one hand, the only element of order $2$ in $\slg(2,K)$ is $-I$, which projects to the identity in $\psl(2,K)$, so the element in $\psl(2,K)$ of order $4$ necessitated by the inclusion must lift to an element of order $8$ in $\slg(2,K)$, implying $\sqrt{2} = z_8\in K$. On the other hand, when $\sqrt{2}\in K$, one can explicitly construct a subgroup $S_4<\psl(2,\mathbb{Q}_2(\sqrt{2}))\leq\psl(2,K)$, as described explicitly in Proposition \ref{prop:exceptional_subgroups}.

    (c) $A_5< \pgl(2,K)$ requires $\sqrt{5}\in K$, so $z_5 = (-1+\sqrt{5})/2\in K$. On the other hand, if $\sqrt{5}\in K$, then $\lvert K/\mathbb{Q}_2\rvert = \lvert K/\mathbb{Q}_2(\sqrt{5})\rvert\cdot \lvert \mathbb{Q}_2(\sqrt{5})/\mathbb{Q}_2\rvert$ is even, so $y^2 + z^2 = -1$ has a solution by (a). This implies $A_5< \pgl(2,K)$ by \cite{serre1972proprietes}.
\end{proof}
\subsection{Cyclic and Dicyclic Subgroups of \texorpdfstring{$\slg(2,K)$}{SL(2,K)}}
A first corollary regarding cyclic subgroups in $\slg(2,K)$ slightly generalizes \cite[Proposition 4.11]{lubotzky1999lattices}:
\begin{cor}[cf. {\cite[Proposition 4.11]{lubotzky1999lattices}}]\label{cor:cyclic}
    Any two cyclic subgroups $A$, $A'< G_x$ of order $n$ with $4\nmid n$ are conjugate in $G_x$.
\end{cor}
\begin{proof}
    This follows from \cite[Proposition 4.11]{lubotzky1999lattices} if $n$ is odd. If $n=2m$ and $m$ is odd, suppose $B<A$, $B'<A'$ are subgroups of order $m$. By \cite[Proposition 4.11]{lubotzky1999lattices}, $B$ and $B'$ are conjugate in $G_x$. Since $A = \langle B,-I\rangle$ and $A' = \langle B',-I\rangle$ and $-I\in Z(G)$, the same $G_x$-conjugation identifies $A$ and $A'$.
\end{proof}

\begin{prop}[cf. {\cite[Proposition 4.6]{lubotzky1999lattices}}]\label{prop:dicyclic_subgroups}
    Let $g\in \slg(2,K)$, $\ord(g) = 2^k$, $k\geq 2$, and $x\in\mathcal{V}(\mathfrak{C}_g)$.
    \begin{itemize}
        \item If $\mathfrak{C}_g$ is a bi-infinite geodesic, there exists a subgroup $A\cong C_{(q-1)2^k}.2< \slg(2,K)_x$ containing $g$. 
        \item If $\mathfrak{C}_g$ is a single vertex, there exists a subgroup $A\cong C_{(q+1)2^k}.2< \slg(2,K)_x$ containing $g$. 
        \item Assume $\mathfrak{C}_g$ is a single edge. If $\lvert K/\mathbb{Q}_2\rvert$ is even, there exists a subgroup $A\cong C_{2^k}.2< \slg(2,K)_x$ containing $g$. If $\lvert K/\mathbb{Q}_2\rvert$ is odd, there is no subgroup $A\leq \slg(2,K)$ with $g\in A$ and $A\cong C_{2^k}.2$.
    \end{itemize}
\end{prop}
\begin{proof}
    (1) By assumption, $K$ contains the root of unity $\xi = \xi_{(q-1)2^k}$. The subgroup can be taken as
    \[
    A = \left\langle \begin{psmallmatrix}
        \xi & 0\\ 0 & \xi^{-1}
    \end{psmallmatrix}, \begin{psmallmatrix}
        0 & -1 \\ 1 & 0
    \end{psmallmatrix}\right\rangle \cong C_{(q-1)2^k}.2.
    \]

    (2) According to the proof of \cite[Theorem 4.2]{lubotzky1999lattices}, it holds that
    \[
    \alpha = \frac{\xi_{q+1}\xi_{2^k}^{-1} - \xi_{q+1}^{-1}\xi_{2^k}}{\xi_{q+1} - \xi_{q+1}^{-1}},\ \beta = \frac{\xi_{2^k} - \xi_{2^k}^{-1}}{\xi_{q+1} - \xi_{q+1}^{-1}}\in K.
    \]
    Therefore,
    \[
    g = \begin{pmatrix}
        \alpha & -\beta \\ \beta & \alpha + \beta z_{q+1}
    \end{pmatrix},\ h = \begin{pmatrix}
        0 & -1 \\ 1 & z_{q+1}
    \end{pmatrix}
    \]
    commute and have orders $2^k$ and $q+1$, meaning $\langle g,h\rangle\cong C_{(q+1)2^k}$.

    Finding an additional generator $\sigma\in\slg(2,K)$ requires
    \[
    \sigma = \begin{pmatrix}
        x & x z_{q+1}+y \\ y & -x
    \end{pmatrix}\in \slg(2,K)\iff x^2 + y^2 + xyz_{q+1} = \lVert x + y\xi_{q+1}\rVert_{K(\xi_{q+1})/K} = -1.
    \]
    Since $K(\xi_{q+1}) = K(\sqrt{-1})$ and $\sigma(\sqrt{-1}) = -\sqrt{-1}$, this is equivalent to finding $x,y\in K$ with $\lVert x + y\sqrt{-1}\rVert_{K(\sqrt{-1})/K} = x^2 + y^2 = -1$. As $K(\xi_{2^k})/K$ does not ramify, $\frac{1}{v_2(\pi_{K})} = \frac{1}{v_2(\pi_{K(\xi_{2^k})})}$, and this common value is even. Therefore $\lvert K/\mathbb{Q}_2\rvert$ is even, and by Lemma \ref{lem:square_sum}, $x^2 + y^2 = -1$ has solutions in $K$. This confirms the existence of $\langle g,h,\sigma\rangle\cong C_{(q+1)2^k}.2<\slg(2,K)_x$.

    (3) Let
    \[
    g = \begin{pmatrix}
        1 & z_{2^k}-2 \\ 1 & z_{2^k}-1
    \end{pmatrix},
    \]
    then $\langle g,\sigma\rangle\cong C_{2^k}.2$ if and only if $\sigma$ takes the form
    \[
    \sigma = \begin{pmatrix}
        x + y & (z-2)x \\ y & -x-y
    \end{pmatrix}\in\slg(2,K),
    \]
    so the existence of $\sigma$ is equivalent to the existence of a solution to
    \[
    -1 = (x+y)^2 - (2-z_{2^k})xy = x^2 + z_{2^k}xy + y^2 = \lVert x + y\xi_{2^k}\rVert_{K(\xi_{2^k})/K}.
    \]
    Again, as $z_{2^k}\in K$, $K(\sqrt{-1}) = K(\xi_{2^k})$ and the existence of $s$ is equivalent to the existence of $x,y\in K$ with $x^2 + y^2 = -1$, which is possible if and only if $\lvert K/\mathbb{Q}_2\rvert$ is even.

    In all cases, since $\langle g\rangle$ is normal, the group $A$ fixes at least one vertex of $\mathfrak{C}_g$.
\end{proof}
% $f(K/\mathbb{Q}_2)>1$ and
\begin{prop}[cf. {\cite[Proposition 4.10]{lubotzky1999lattices}}]\label{prop:dicyclic_fixed}
    Assume $A\leq \slg(2,K)$ is a subgroup isomorphic to $C_{2^k}.2$.
    \begin{itemize}
        \item There exists a unique $\gl(2,K)$-conjugacy class of subgroups isomorphic to $A$.
        \item If $r_k$ is even, $T_A$ is an $r_k/2$-neighborhood of a vertex. If $r_k$ is odd, $T_A$ is an $(r_k-1)/2$-neighborhood of an edge. For the generator $g$ of the $C_{2^k}$-subgroup in $A$, $\mathfrak{C}_g\cap T_A$ consists of a single vertex if $\mathfrak{C}_g$ is either an infinite path or a vertex, and $\mathfrak{C}_g\cap T_A = \mathfrak{C}_g$ if it is a single edge.
    \end{itemize}
\end{prop}
\begin{proof}
    (1) This follows by the same argument as in the proof of \cite[Proposition 4.10 (1)]{lubotzky1999lattices}.

    (2) (a) If $\xi_{2^k}\in K$, as in Proposition \ref{prop:dicyclic_subgroups}, we may take
    \[
    g = \begin{pmatrix}
        \xi_{2^k} & 0\\ 0 & \xi_{2^k}^{-1}
    \end{pmatrix},\ \sigma = \begin{pmatrix}
        0 & -1 \\ 1 & 0
    \end{pmatrix},\ A = \langle g,\sigma\rangle.
    \]
    The axis of $g$ joins $0$ and $\infty\in\mathbf{P}^1(K)$, and that of $\sigma$ joins $\pm \sqrt{-1} \in\mathbf{P}^1(K)$. By Theorem \ref{thm:regular_width}, we have the fixed trees $T_g = N(\mathfrak{C}_g, r_k)$, $T_\sigma = N(\mathfrak{C}_\sigma,r_2)$, where the cores $\mathfrak{C}_g,\mathfrak{C}_\sigma$ are the axes, with distance between them given by
    \[
    d(\mathfrak{C}_g,\mathfrak{C}_\sigma) = \frac{v(2)}{v(\pi_K)} = r_2.
    \]
    Since $k\geq 2$ and $r_k = r_2/2^{k-2}\leq r_2$, we have
    \[
    T_A = T_g\cap T_\sigma = N(\mathfrak{C}_g, r_k)\cap N(\mathfrak{C}_\sigma,r_2) = N(m,r_k/2),
    \]
    where $m$ is the midpoint of a path of length $r_k$ contained in the path of length $r_2$ connecting $\mathfrak{C}_g$ and $\mathfrak{C}_\sigma$. It is either the $r_k/2$-neighborhood of a vertex or the $(r_k-1)/2$-neighborhood of an edge, depending on the parity of $r_k$.

    (b.i) If $\xi_{2^k}\notin K$ and $K(\xi_{2^k})/K$ is unramified, we may take
    \[
    g = \begin{pmatrix}
        \alpha & -\beta \\ \beta & \alpha + \beta z_{q+1}
    \end{pmatrix},\ \sigma = \begin{pmatrix}
        x & x z_{q+1}+y \\ y & -x
    \end{pmatrix},
    \]
    where
    \[
    \alpha = \frac{\xi_{q+1}\xi_{2^k}^{-1} - \xi_{q+1}^{-1}\xi_{2^k}}{\xi_{q+1} - \xi_{q+1}^{-1}},\ \beta = \frac{\xi_{2^k} - \xi_{2^k}^{-1}}{\xi_{q+1} - \xi_{q+1}^{-1}}\in K,
    \]
    and $x,y\in K$, $x^2+z_{q+1}xy + y^2 = -1$. Solutions $x,y$ correspond to solutions of $x'^2 + y'^2 = -1$ by
    \[
    x = x' - \frac{z_{q+1}}{\gamma}y',\ y = \frac{2}{\gamma}y',
    \]
    where
    \[
    \gamma = \frac{\xi_{q+1} - \xi_{q+1}^{-1}}{\sqrt{-1}}\in \mathcal{O}_K^\times,
    \]
    % For $k = 2$, $z = 0$ and $\gamma = 2$, while for $k\geq 3$, $v(z) = 1$ and $v(\gamma) = \frac{1}{2^{k-2}}$. In both cases,
    so both $\frac{z}{\gamma}$ and $\frac{2}{\gamma}$ are in the ring of integers $\mathcal{O}_K$. By Lemma \ref{lem:square_sum}, $x'^2 + y'^2 = -1$ has integral solutions, so the corresponding $x,y$ with $x^2 + z_{q+1}xy + y^2 = -1$ are also integral.

    Since $v(\beta) = \frac{1}{2^{k-2}} = r_kv(\pi)$, $g$ fixes a subpath of $(0,\infty)$ of length $2r_k$ centered at the origin; indeed the origin is the core of $g$. On the other hand, $v(z_{q+1}) = 0$, and since at least one of $x$ and $y$ has zero valuation, either $y$ or $xz_{q+1}+y$ has zero valuation. This shows that the origin is in $\partial T_{\sigma}$. By a similar argument, $T_A = T_g\cap T_\sigma$ is either the $r_k/2$-neighborhood of a vertex or the $(r_k-1)/2$-neighborhood of an edge.

    (b.ii) Since $K(\xi_{2^k})=K(\sqrt{-1})$ is ramified over $K$, both $\mathfrak{C}_g$ and $\mathfrak{C}_\sigma$ are single edges. Thus $T_A=T_g\cap T_\sigma$ is the intersection of two edge-neighborhoods. On the vertex set, each edge-neighborhood is a ball centered at the midpoint of the edge, and the nonempty intersection of two balls in a tree is again a ball. Consequently, $T_A$ is a neighborhood of either a vertex or an edge. It therefore suffices to show that its diameter is $r_k$.
    
    On the one hand, the diameter of $T_A$ is at least $r_k$: setting
    \[
    g = \begin{pmatrix}
        1 & z_{2^k}-2 \\ 1 & z_{2^k}-1
    \end{pmatrix},\ \sigma = \begin{pmatrix}
        x + y & (z-2)x \\ y & -x-y
    \end{pmatrix},
    \]
    where $x,y\in \mathcal{O}_K$, we see that both matrices lie in the following vertex stabilizers:
    \[
    \left\{\begin{pmatrix}
        a & b\\ c & d
    \end{pmatrix}\in\slg(2,K)\mid a,d\in\mathcal{O}_K,\ b\in \pi^i \mathcal{O}_K,\ c\in \pi^{-i}\mathcal{O}_K\right\},\ i=0,\dots, v(z-2) = r_k,
    \]
    that is, $T_g\cap T_\sigma = T_A$ contains a path of length $r_k$.

    On the other hand, the diameter of $T_A$ cannot be larger. If $r_k$ is odd, then by Theorem \ref{thm:width_formula}, the diameter of $T_g$ is already $r_k$. If $r_k$ is even, up to $\gl(2,K)$-conjugation, we may assume that $T_g$ contains a subpath of $(0,\infty)$ of length $2d_k+1$, where $\frac{r_k}{2}\leq d_k\leq r_k-1$, corresponding to the interval $[-(2d_k+1),0]$ in the standard apartment; we suppose to the contrary that $T_\sigma\cap T_g$ contains a subpath of $(0,\infty)$ of length $l>r_k$. It follows that
    \[
    g = \begin{pmatrix}
        \alpha & \epsilon \\ \epsilon^{-1}(\alpha(z-\alpha)-1) & z-\alpha
    \end{pmatrix},\ \epsilon\in\mathcal{O}_K^\times,
    \]
    with $v_K(\alpha(z-\alpha)-1) = 2d_k+1$, and the relation $\sigma g = g^{-1}\sigma$ implies
    \[
    \sigma = \begin{pmatrix}
        \frac{\epsilon^{-1}(\alpha(z-\alpha)-1)x + \epsilon y}{z-2\alpha} & x\\ y & - \frac{\epsilon^{-1}(\alpha(z-\alpha)-1)x + \epsilon y}{z-2\alpha}
    \end{pmatrix}.
    \]
    Since $\sigma$ fixes this common subpath of length $l$, we have
    \[
    v_K(y)\geq l,\ v_K(x)\geq l-(2d_k+1),\text{ and }v_K(xy)\geq l> r_k.
    \]
    Consequently, 
    \[
    v_K\bigl(
    \epsilon^{-1}(\alpha(z-\alpha)-1)x+\epsilon y
    \bigr)\geq l>r_k.
    \]
    However, from the proof of Theorem \ref{thm:width_formula} we see that $v(\alpha - 1)>0$, so $v_K(2(1-\alpha))> v_K(2)\geq r_k$. Together with $v_K(z-2) = r_k$, this implies
    \[
    v_K(z-2\alpha) = r_k.
    \]
    Therefore,
    \[
    v_K\left(- \frac{\epsilon^{-1}(\alpha(z-\alpha)-1)x + \epsilon y}{z-2\alpha}\right)>0,
    \]
    so $v(\det(\sigma))>0$, contradicting the assumption that $\sigma \in \slg(2,K)$.
    
\end{proof}
Unlike in the odd-prime cases, the group $N_{\gl(2,K)}(A)$ need not act transitively on $\partial T_A$ for $p=2$ (see Appendix \ref{append:b}). As a result, an analogue of \cite[Proposition 4.10, (2)]{lubotzky1999lattices} no longer holds: for $x\in \mathcal{V}(X)$, there are at least $\lvert \partial T_A\rvert/2$ (or $\lvert \partial T_A\rvert/6$ when $k=2$) $\gl(2,K)_x$-conjugacy classes of groups $A^*\cong C_{2^k}.2$ with $x\in \partial T_{A^*}$ such that $A^*\cong A$. Consequently, the argument of \cite[Lemma 4.7(3)]{lubotzky1999lattices}, which plays a key role in their lattice constructions, does not carry over directly. Given a subgroup $A\cong C_{2^k}.2$ of $\slg(2,K)_x$ with $x\in\partial T_A$, the existence of a finite subgroup $B<\slg(2,K)_x$ containing a subgroup isomorphic to $A$ does not in general guarantee the existence of a subgroup $B^*<\slg(2,K)_x$ containing $A$ with $B^*\cong B$. This leads to the following question:
\begin{itemize}
    \item[(*)] Let $A\cong C_{2^k}.2$ be a subgroup of $\slg(2,K)$. For each vertex $x\in\partial T_A$, which isomorphism types of finite proper supergroups $B>A$ can occur with $B<\slg(2,K)_x$?
\end{itemize}
By Theorem \ref{thm:finite_subgroup}, any such $B$ is isomorphic to $C_{s\cdot 2^l}.2$, where $s$ is odd, $l\geq k$, and $s\mid q-1$ or $s\mid q+1$; to $2.A_4$ or $2.A_5$ when $k=2$; or to $2.S_4$ when $k\leq 3$. The following proposition handles the case when $B$ is dicyclic.
\begin{prop}\label{prop:supergroup}
    Fix a subgroup $A<\slg(2,K)$ with presentation $A = \langle g,\sigma\mid g^{2^{k-1}} = \sigma^2 = -I,\ \sigma g\sigma^{-1} = g^{-1}\rangle\cong C_{2^k}.2$. In the following, consider only integers $l\geq k$ and odd integers $s\geq 1$ such that $C_{s\cdot 2^l}.2$ is isomorphic to a subgroup of $\slg(2,K)$, and $(s,l)\neq (1,k)$.
    
    If $k\geq 3$, then a vertex $x\in\partial T_A$ allows a supergroup $B>A$, $B<\slg(2,K)_x$ and $B \cong C_{2^l}.2$ if and only if $d(x,\mathfrak{C}_g)\leq d_l$, where $d_l$ is the width defined in Theorem \ref{thm:regular_width}. For $s>1$, a vertex $x\in\partial T_A$ allows a supergroup $B>A$, $B<\slg(2,K)_x$ and $B \cong C_{s\cdot 2^l}.2$ if and only if $d(x,\mathfrak{C}_g) = 0$.

    If $k = 2$, a similar conclusion holds with $\mathfrak{C}_g$ being replaced by $\mathfrak{C}_g\cup \mathfrak{C}_\sigma \cup \mathfrak{C}_{g\sigma}$.
\end{prop}
\begin{proof}
    For $k\geq 3$, there is at most one subgroup $B\cong C_{s\cdot 2^l}.2$ containing $A = \langle g,\sigma\rangle\cong C_{2^k}.2$. Indeed, the cyclic subgroup of index two in $B$ must contain $g$, so it is the unique cyclic subgroup of order $s\cdot 2^l$ containing $g$. Suppose $h$ is a generator. Then $B$ is generated by $\sigma$ and $h$. If $G$ contains subgroups isomorphic to $C_{s\cdot 2^l}.2$, then the existence of such a supergroup $B>A$ is guaranteed by the uniqueness of the conjugacy class in Proposition \ref{prop:dicyclic_fixed}. Since $T_B = T_h\cap T_\sigma$ and $\partial T_A\subset T_\sigma$, a vertex $x\in\partial T_A$ allows $B<G_x$, or equivalently $x\in T_B$, if and only if $x\in T_h$. 

    If $s=1$, then $\mathfrak{C}_h=\mathfrak{C}_g$ and $T_h=N(\mathfrak{C}_g,d_l)$, so $B<G_x$ if and only if $d(x,\mathfrak{C}_g)\leq d_l$. If $s>1$, then $h^{2^l}$ fixes either a vertex or a bi-infinite path, and Lemma \ref{lem:odd_fixed} and Theorem \ref{thm:regular_width} imply its fixed subtree is the same as $\mathfrak{C}_g$. Therefore, $T_h = \mathfrak{C}_g$, and $B<G_x$ if and only if $d(x,\mathfrak{C}_g) = 0$.

    For $k = 2$, the supergroup $B>A\cong Q_8$ can arise in three ways, according to which of $g,\sigma$ and $g\sigma\in A$ lies in the cyclic subgroup $C_{s\cdot 2^l}<B$. Considering all three possibilities, the core $\mathfrak{C}_g$ is replaced by the union $\mathfrak{C}_g\cup \mathfrak{C}_\sigma \cup \mathfrak{C}_{g\sigma}$ in the corresponding conclusion.
\end{proof}
To discuss the cases when $B$ is binary polyhedral, we need the following lemma describing the structure of $T_{Q_8}$ and $T_{Q_{16}}$.
\begin{lem}\label{lem:exceptional_tree_structure}
    Fix the presentation $A = \langle g,\sigma\mid g^{2^{k-1}} = \sigma^2 = -I,\ \sigma g\sigma^{-1} = g^{-1}\rangle\cong C_{2^k}.2$.
    \begin{itemize}
        \item If $k=2$ and a subgroup of $\slg(2,K)$ is isomorphic to $2.A_4$, then there exists a unique element $q\in \slg(2,K)$, $\ord(q) = 3$, such that $qgq^{-1} = \sigma$ and $q^{-1}gq = g\sigma$. Moreover, $T_q$ is either the central vertex of $T_A$ or a bi-infinite geodesic passing through the center of $T_A$.
        \item If $k=3$ and a subgroup of $\slg(2,K)$ is isomorphic to $2.S_4$, then there exist unique elements $q_i\in \slg(2,K)$, $\ord(q_i) = 3$, $i=1,2$, such that $q_ig^2q_i^{-1} = g^{i-1}\sigma$ and $q_i^{-1}g^2q_i = g^{i+1}\sigma$. Moreover, for $i=1,2$, $T_{q_i}$ meets $T_A$ exactly at a vertex in $\partial T_A$, and the distance between these vertices is $r_3$.
    \end{itemize}
\end{lem}
\begin{proof}
    (1) If $A\cong Q_8$ and $\slg(2,K)$ admits a subgroup $B\cong 2.A_4$, then $B$ contains a subgroup $A^*\cong A$. Since $A$ has a unique $\gl(2,K)$-conjugacy class (Proposition \ref{prop:dicyclic_fixed}), the corresponding adjoint homomorphism sends $B$ to a supergroup $B^*>A$. Then $B^*$ contains an element $q$ satisfying the relations above. Since the centralizer $C_G(A)$ contains only scalars, $q$ is unique up to $\pm 1$, and exactly one of these two choices has order $3$. According to Lemma \ref{lem:odd_fixed}, $q$ fixes either a bi-infinite path or a single vertex, depending on the parity of $f(K/\mathbb{Q}_2)$. Moreover, $q\in N_G(A)$ implies $q$ preserves $T_A$, so $T_q$ is either the center of $T_A$, or a path through the center.

    (2) If $A\cong Q_{16}$, then it contains exactly two subgroups $Q,Q'<A$ isomorphic to $Q_8$, namely $\langle g^2,\sigma\rangle$ and $\langle g^2,g\sigma\rangle$. By (1), $Q$ and $Q'$ correspond to elements $q,q'$ of order $3$ with the prescribed relationship to $A$. In addition, since $\langle Q,Q'\rangle = A$, we have $T_Q\cap T_{Q'} = T_A$. Proposition \ref{prop:dicyclic_fixed} implies that both $T_Q$ and $T_{Q'}$ are balls of diameter $r_2$, while $T_A$ is a ball of diameter $r_3$. Since $r_3 = r_2/2$, the centers of $T_Q$ and $T_{Q'}$ are at distance $r_3$ from each other, and the midpoint is the center of $T_A$. Therefore, both centers lie on $\partial T_A$. For the case where $T_q$ is an infinite path, since $r_2 = 2r_3$ is even, the three cores of $g^2, \sigma$ and $g^2\sigma\in Q$ are distinct, and conjugation by $q$ cyclically permutes them. The direction from the center of $T_Q$ toward the center of $T_A$, which is the direction toward the core of $g$, is therefore not fixed by $q$. Hence, $T_q$ and $T_A$ meet at exactly a single vertex.
\end{proof}
\begin{prop}\label{prop:supergroup_exceptional}
    If $k=3$, then $x\in\partial T_A$ allows a supergroup $A<B<\slg(2,K)_x$ with $B\cong 2.S_4$ if and only if $x\in T_{q_i}\cap \partial T_A$, $i=1,2$.
    
    If $k=2$, then $x\in\partial T_A$ allows a supergroup $A<B<\slg(2,K)_x$ with $B\cong 2.A_4$ or $B\cong 2.A_5$ if and only if $2\mid f$ and $x\in T_q\cap \partial T_A$. 
    
    The vertex $x\in\partial T_A$ allows a supergroup $B<\slg(2,K)_x$ isomorphic to $2.S_4$ if and only if $x \in \partial T_A\cap T_{q_i}$, $i=1,2,3$, where $q_i$ denote the elements different from $q$ in the three supergroups of $A$ isomorphic to $Q_{16}$ described in Lemma \ref{lem:exceptional_tree_structure}. $A$ is not normal in these supergroups $B$.
\end{prop}
\begin{proof}
    (1) Suppose $k=3$ and $A\cong Q_{16}$. Then $2.S_4$ is the only binary polyhedral group that may contain $A$. If $G$ admits a subgroup isomorphic to $2.S_4$, there is a supergroup $B>A$ with $B\cong 2.S_4$. Now, $B$ contains a normal subgroup $Q\cong Q_8$, and since $A$ is a Sylow $2$-group of $B$, we have $Q\leq A$. Therefore, $Q$ is one of the two $Q_8$-subgroups, and the proof of Lemma \ref{lem:exceptional_tree_structure} implies that either $B = \langle A,q_1\rangle$ or $B = \langle A,q_2\rangle$. In either case, $T_B = T_A\cap T_{q_i}$ is a vertex in $\partial T_A$. Both possibilities are allowed, as $I+g\in\gl(2,K)$ normalizes $A$, interchanges its two $Q_8$-subgroups, and thus interchanges the two $2.S_4$-supergroups.
    
    (2) Suppose $k=2$ and $A\cong Q_8$. If $G$ admits a subgroup isomorphic to $2.A_4$, we have shown in Lemma \ref{lem:exceptional_tree_structure} that there is a supergroup $B>A$ in the isomorphism class, and the uniqueness of $q$ implies the uniqueness of $B = \langle q,A\rangle$. Since $T_B = T_A\cap T_q$, $x\in T_B\cap \partial T_A$ if and only if $T_q$ is an infinite path, i.e., $2\mid f$, and $x$ is one of the two intersection vertices in $T_q\cap \partial T_A$.

    Similarly, if $G$ admits a subgroup isomorphic to $2.A_5$, the uniqueness of $Q_8$-conjugacy classes provides a supergroup $B>A$ with $B\cong 2.A_5$. There is an intermediate group $A<2.A_4<B$, and since $2\mid f$, we have shown that $T_{2.A_4}$ is a diameter of $T_A$, and $T_B$ is a subset of it. If $B$ fixes a vertex $y\notin \partial T_A$, $g\in A$ fixes the star of $y$, and so does the conjugation of $g$ by an order-$5$ element. Since these elements generate the entire $B\cong 2.A_5$, $\st(y)\subset T_B$, a contradiction. Since $T_B$ is non-empty, it is one of the vertices of $\partial T_{2.A_4} = \partial T_A\cap T_q$. Both possibilities are allowed, as $I+g\in\gl(2,K)$ normalizes $A$ and interchanges the cosets $qQ_8$ and $q^{-1}Q_8$, and thus interchanges the two endpoints of $T_A\cap T_q$.

    Now, if $A<B$ with $B\cong 2.S_4$, there is an intermediate group $A<Q_{16}<B$. There are three supergroups $Q_i>A$, $i=1,2,3$ isomorphic to $Q_{16}$, and by (1), each $Q_i$ corresponds to two supergroups isomorphic to $2.S_4$. By Lemma \ref{lem:exceptional_tree_structure}, one of them contains an order-$3$ element $q$ such that $\langle A, q\rangle\cong 2.A_4$, so this element is the same for all three $Q_i$. In this supergroup, $A$ is normal, and its fixed subtree is $T_q\cap T_{Q_i}$, which consists of the center of $T_A$ and hence contains no boundary vertex of $T_A$. 
    
    Suppose the other supergroup is $B_i = \langle Q_i,q_i\rangle$, with $\ord(q_i)=3$. Then the normal subgroup $A_i\triangleleft B_i$, $A_i\cong Q_8$ is different from $A$. By (1), $B_i$ fixes the center of the fixed subtree of $A_i$, and by Lemma \ref{lem:exceptional_tree_structure}, the distance between the centers of the fixed subtrees of $A$ and $A_i$ is $r_3 = r_2/2$. Therefore, $B_i$ fixes exactly one vertex on $\partial T_A$. These vertices lie in $T_{Q_i}$ but are different from the center of $T_A$. Since $T_q\cap T_{Q_i}$ consists only of the center of $T_A$, none of them lies in $T_q$. They are cyclically permuted by $q$ and are therefore distinct.
\end{proof}
\begin{exm}
    First, suppose $K=\mathbb{Q}_2(z_8,z_{20})$, so $q=4$, $e(K/\mathbb{Q}_2)=r_2=4$, and $K(i)/K$ is quadratic and unramified. For each subgroup $A<\slg(2,K)$ with $A\cong Q_8$, the fixed subtree $T_A$ is the $2$-neighborhood of a vertex $v$. The three cores corresponding to the cyclic subgroups of order $4$ in $A$ are vertices $u_i\in\partial T_A$, $i=1,2,3$. The paths $[v,u_i]$ are diameters of the fixed subtrees $T_{Q_i}$ of the three supergroups $Q_i>A$ with $Q_i\cong Q_{16}$. Each $Q_i$ contains another subgroup $A_i\cong Q_8$, $A_i\neq A$, whose fixed subtree has center $v_i\in\partial T_A$. The vertices $v,u_i,v_i$ are pairwise at distance $2$.

    The fixed subtree of the unique $2.A_4$ supergroup of $A$ is a diameter $P$ of $T_A$. The diameter endpoints $w,w'$ are the unique fixed vertices of the two supergroups isomorphic to $2.A_5$ containing this $2.A_4$, respectively. The unique $2.S_4$ supergroup in which $A$ is normal fixes exactly $v$, while the three $2.S_4$ supergroups in which $A$ is not normal fix exactly $v_1,v_2,v_3$, respectively.

    Next, suppose $K = \mathbb{Q}_2(z_8,z_5)$; then $e(K/\mathbb{Q}_2) = 2$, and $K(i)/K$ is now ramified. In this case, $T_A$ is a ball of radius $1$, and the cores of order-$4$ elements become the edges connecting the normal $2.S_4$ fixed vertex $v$ and the non-normal $2.S_4$ fixed vertices $v_i$, $i=1,2,3$. The other two edges form the $2.A_4$ fixed path. 
    
    The fixed subtrees $T_A$ in both cases, together with these fixed subtrees and cores, are illustrated below.
    \begin{center}
        \begin{minipage}{0.48\textwidth}
            \centering
            \includegraphics[trim={130 140 120 120}, clip, scale=0.7]{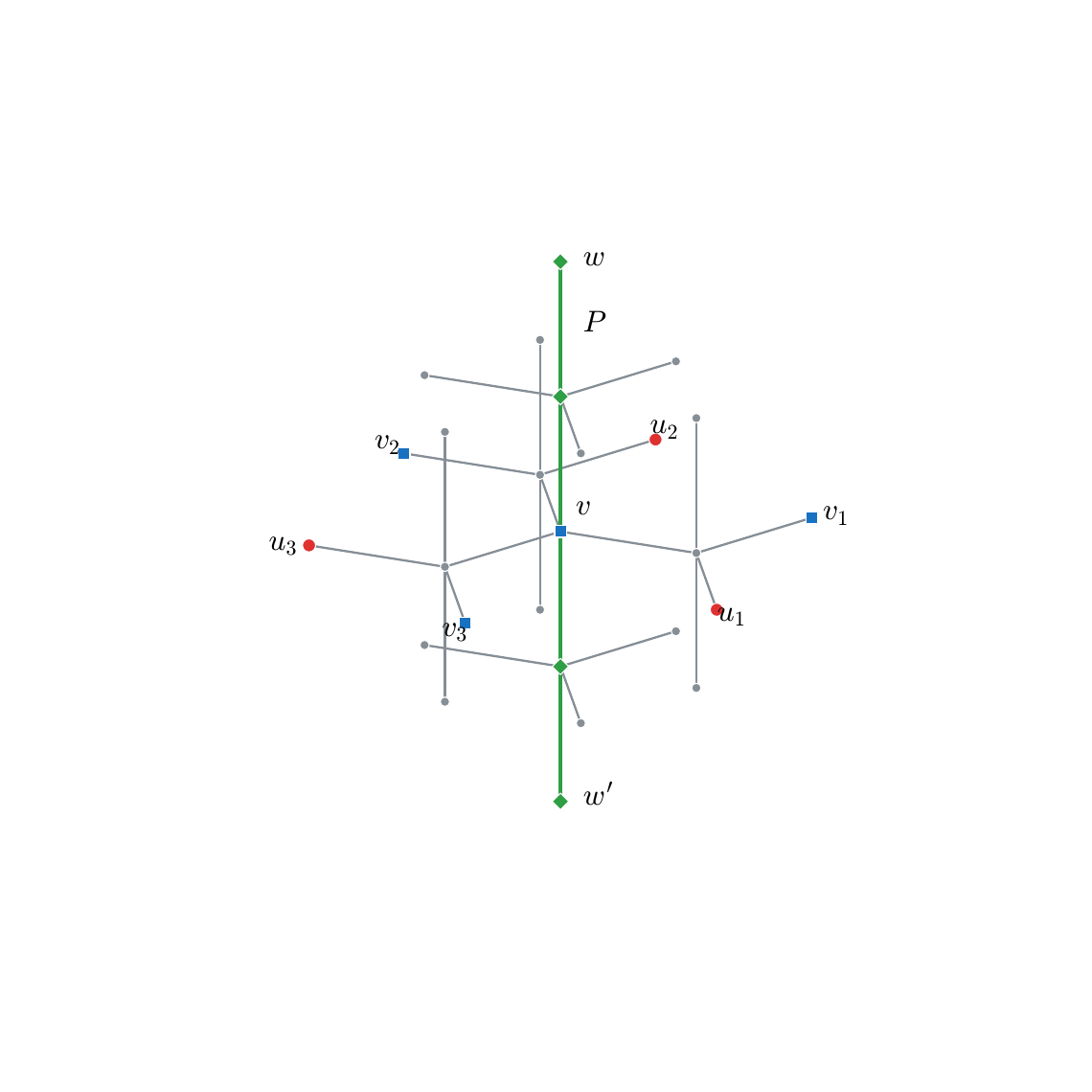}
        \end{minipage}
        \hfill
        \begin{minipage}{0.48\textwidth}
            \centering
            \includegraphics[trim={130 140 120 120}, clip, scale=0.7]{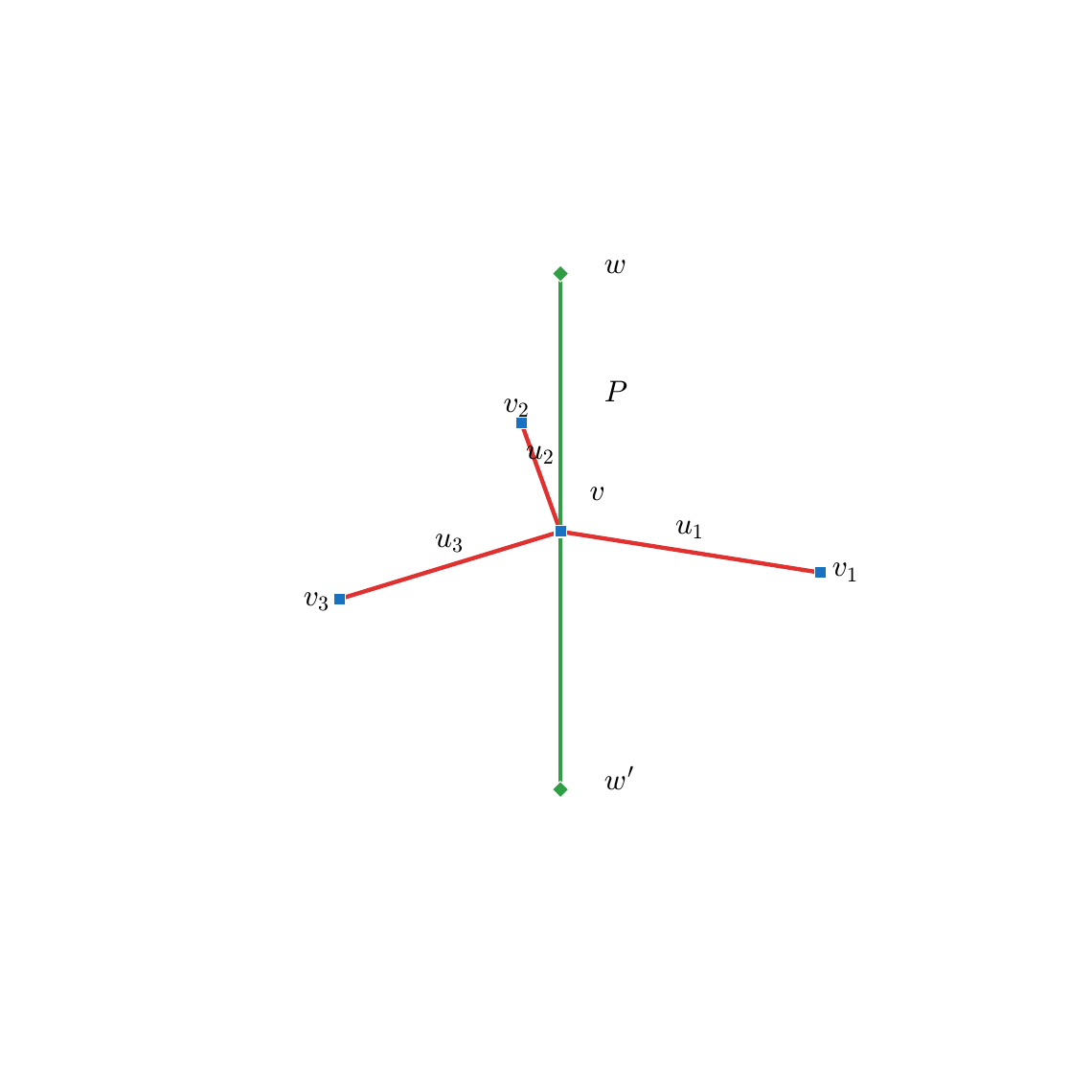}
        \end{minipage}
        % left bottom right top
    \end{center}
\end{exm}
\subsection{Binary Polyhedral Subgroups}
\begin{prop}\label{prop:exceptional_subgroups}
    Let $g\in \slg(2,K)$ satisfy $\ord(g) = 4$, and suppose that $\lvert K/\mathbb{Q}_2\rvert$ is even.
    \begin{itemize}
        \item There exists $x\in\mathcal{V}(T_g)$ such that there is a subgroup $A\cong 2.A_4< \slg(2,K)_x$ containing $g$.
        \item If $\sqrt{2}\in K$, there also exists a subgroup $A\cong 2.S_4< \slg(2,K)_x$ containing $g$.
    \end{itemize}
\end{prop}
\begin{proof}
    Lemma \ref{lem:square_sum} implies the existence of $x,y\in K$ such that $x^2 + y^2 = -1$. Define the following matrices:
    \[
    g = \begin{pmatrix}
        0 & 1 \\ -1 & 0
    \end{pmatrix},\ q = \frac{1}{2}\begin{pmatrix}
        1 - x + y & 1 + x + y \\ -1 + x + y & 1 + x - y
    \end{pmatrix},
    \]
    so
    \[
    \det(g) = 1,\ \det(q) = \frac{1}{4}((1-(x-y)^2) - ((x+y)^2 - 1)) = \frac{1}{4}(2 - 2(x^2 + y^2)) = 1,
    \]
    and
    \[
    \tr(g) = 0,\ \tr(q) = 1,
    \]
    so $\ord(g) = 4$ and $\ord(q) = 6$. At the same time,
    \[
    gq = \frac{1}{2}\begin{pmatrix}
        -1+x+y & 1+x-y \\ -1+x-y & -1-x-y
    \end{pmatrix},\ \tr(gq) = -1,
    \]
    so $\ord(gq) = 3$. Therefore, the group generated by $\overline{g},\overline{q}\in \psl(2,K)$ is
    \[
    \langle \overline{g},\overline{q}\mid \overline{g}^2, \overline{q}^3, (\overline{g}\overline{q})^3\rangle\cong A_4.
    \] 
    If $\sqrt{2}\in K$, we can take
    \[
    x = 2+\sqrt{\frac{-6}{5}},\ y = 1-2\sqrt{\frac{-6}{5}}\in K,
    \]
    as $-\frac{3}{5}\equiv 1(\bmod\ 8)$. 

    In addition, we define
    \[
    p = \frac{\sqrt{2}}{2}\begin{pmatrix}
        1 & 1 \\ -1 & 1
    \end{pmatrix},
    \]
    so $p^2 = g$ and $\tr(p) = \sqrt{2}$, and $\ord(p) = 8$. At the same time,
    \[
    pq = \frac{\sqrt{2}}{2}\begin{pmatrix}
        y & x + 1 \\ x - 1 & -y
    \end{pmatrix},\ \tr(pq) = 0,
    \]
    so $\ord(pq) = 4$. Therefore, the group generated by $\overline{p},\overline{q}\in \psl(2,\mathbb{Q}_2(\sqrt{2}))$ is
    \[
    \langle \overline{p},\overline{q}\mid \overline{p}^4, \overline{q}^3, (\overline{p}\overline{q})^2\rangle\cong S_4.
    \]
    These groups lift to $2.A_4$ and $2.S_4$ in $\slg(2,K)$.
\end{proof}
\begin{prop}\label{prop:x_to_core}
    In the setting of Proposition \ref{prop:exceptional_subgroups}, suppose $x$ is taken such that $d(x,\mathfrak{C}_g)$ is minimized.
    \begin{itemize}
        \item If $\mathfrak{C}_g$ is either a bi-infinite geodesic or a single vertex, then $d(x,\mathfrak{C}_g) = r_2/2$.
        \item If $\mathfrak{C}_g$ is an edge and $r_2$ is even, then $d(x,\mathfrak{C}_g) = d_2 - r_2/2$.
        \item If $\mathfrak{C}_g$ is an edge and $r_2$ is odd, then $x\in \mathfrak{C}_g$.
        % the binary polyhedral subgroup $A< \slg(2,K)_x$ fixes $\mathfrak{C}_g$.
    \end{itemize}
    % In addition, the 
\end{prop}
\begin{proof}
    Consider the non-trivial elements in the subgroup $Q\triangleleft A$, $Q\cong Q_8$, namely
    \[
    g = p^2 = \begin{pmatrix}
        0 & 1 \\ -1 & 0
    \end{pmatrix},\ qgq^{-1} = \begin{pmatrix}
        -x & y \\ y & x
    \end{pmatrix},\ q^{-1}gq =\begin{pmatrix}
        y & x \\ x & -y
    \end{pmatrix}.
    \]
    We notice that the intersection of $T_g = N(\mathfrak{C}_g,d_2)$, $T_{qgq^{-1}} = N(\mathfrak{C}_{qgq^{-1}},d_2)$ and $T_{q^{-1}gq} = N(\mathfrak{C}_{q^{-1}gq},d_2)$ coincides with $T_Q$, which is a ball of diameter $r_2$ by Proposition \ref{prop:dicyclic_fixed}. 
    
    Suppose first that $\mathfrak C_g$ is a bi-infinite geodesic or a single vertex. The three cores $\mathfrak{C}_g$, $\mathfrak{C}_{qgq^{-1}}$ and $\mathfrak{C}_{q^{-1}gq}$ are cyclically permuted by $q$, and since $d_2 = r_2$ in these cases, they cannot have a common vertex. By the symmetry of these cores, they are connected by a tripod. Since $A$ permutes these cores, the tripod center $x$ is fixed by $A$, and is the vertex in $T_A$ closest to $\mathfrak C_g$. Using the diameter of $T_Q$ and the width of $T_g$, it is straightforward to see that $d(x,\mathfrak{C}_g) = r_2/2$.
    
    Now we consider the case when $\mathfrak{C}_g$ is an edge. If the three edges coincide, then the diameter of $T_Q=N(\mathfrak{C}_g,d_2)$ is $r_2 = 2d_2+1$, which is odd. Since $Q\triangleleft A$, the group $A$ preserves the central edge $\mathfrak{C}_g$ and therefore fixes it pointwise. Therefore, $x$ is one of the vertices of $\mathfrak{C}_g$.
    
    If the three edges are distinct, then the subtree spanned by them is a tripod. As $Q\triangleleft A$, $A$ permutes the three core edges and hence fixes their tripod center $x$. Thus, it is the vertex in $T_A$ closest to $\mathfrak{C}_g$. Since $T_Q$ is the intersection of $T_g$ with its two images under $q$, $x$ is the center of $T_Q$. Using the diameters of $T_g$ and $T_Q$, it is straightforward to see that $r_2$ is even, and $d(x,\mathfrak{C}_g) = d_2 - r_2/2$.
\end{proof}
Unlike in the odd-$p$ cases, the binary polyhedral group $2.A_4$ may occur as an edge group for minimal covolume lattices of $\slg(2,K)$ with $K/\mathbb{Q}_2$, e.g., Example \ref{exm:3}.
% For example, the construction in Section \ref{sec:construction} gives a lattice
% \[
% \Gamma\cong 2.A_5 *_{2.A_4} 2.A_5,\ \Gamma< \slg(2,\mathbb{Q}_2(\xi_3)),
% \]
% with covolume $\frac{1}{120}$. 
For this reason, we similarly discuss the fixed subtree $T_A$ for $A\cong 2.A_4$ and how supergroups $B>A$ can be attached at vertices $x\in \partial T_A$.
\begin{cor}\label{cor:A4_attach}
    Assume $A\leq \slg(2,K)$ is a subgroup isomorphic to $2.A_4$, and let $Q$ be its unique normal subgroup isomorphic to $Q_8$.
    \begin{itemize}
        \item There exists a unique $\gl(2,K)$-conjugacy class of subgroups isomorphic to $A$.
        \item If $2\mid f$, then $T_A$ is a path of length $r_2$ sharing its center with $T_Q$. If $2\nmid f$, then $T_A$ consists of the central vertex of $T_Q$.
        \item If $2\mid f$, then for each $x\in\partial T_A$, there exists a supergroup $A<B<\slg(2,K)_x$ with $B\cong 2.A_5$, whereas no such supergroup is isomorphic to $2.S_4$.
    \end{itemize}
\end{cor}
\begin{proof}
    The first claim follows from the proof of Lemma \ref{lem:exceptional_tree_structure}, and the second claim follows from the proof of Proposition \ref{prop:supergroup_exceptional}. The third claim also follows from the proof of Proposition \ref{prop:supergroup_exceptional}: for $B\cong 2.A_5$, the two vertices $x\in T_q\cap \partial T_Q$ allowing $B<G_x$ are just the vertices in $\partial T_A$. If $B\cong 2.S_4$, then $Q\triangleleft B$, so the non-normal requirement in Proposition \ref{prop:supergroup_exceptional} implies that $B<G_x$ is not possible for every $x\in \partial T_A\subset\partial T_Q$.
\end{proof}
\subsection{Restrictions on Vertex Stabilizers}

\begin{prop}[cf. {\cite[Proposition 3.4]{lubotzky1999lattices}}]\label{prop:120_bound}
    Let $g\in \slg(2,\mathcal{O}_K)$ be an element of order $2^k$, $k\geq 2$, such that
    \[
    g\notin \ker(\tau:\ \slg(2,\mathcal{O}_K)\to \slg(2,\mathbb{F}_q)).
    \]
    Then for a finite subgroup $S < \slg(2,\mathcal{O}_K)$\footnote{\cite[Proposition 3.4]{lubotzky1999lattices} appears to have a typo here: their statement with the assumption $S < \slg(2,K)$ does not hold; in fact, this assumption is inconsistent with the context and is not needed in this generality later.} containing $g$, one of the following holds:
    \begin{itemize}
        \item The group $C_{\slg(2,\mathcal{O})}(g)\cap S$ is a cyclic group of order $2^{k}$ and has index at most $2$ in $S$.
        \item $k = 2$, $C_{\slg(2,\mathcal{O})}(g)\cap S = \langle g\rangle = C_4$, and $g$ is not contained in the index $2$ cyclic subgroup of $S$ described in Theorem \ref{thm:finite_subgroup}.
        \item $k \in \{2, 3\}$ and $S\cong 2.A_4$, $2.S_4$, or $2.A_5$.
    \end{itemize}
\end{prop}
\begin{proof}
    When $S$ is not isomorphic to $2.A_4$, $2.S_4$, or $2.A_5$, it is contained in a maximal subgroup of type (a) or (b) as described in Theorem \ref{thm:finite_subgroup}.

    Note that the maximal subgroup $\mu_K.2$ or $\mu_{E/K}^*.2$ is a dicyclic group. When $g$ is not contained in the cyclic subgroup $\mu_K$ or $\mu_{E/K}^*$, it holds that $\ord(g) = 4$, and its centralizer in the dicyclic group is the order-four subgroup generated by it, and so is its intersection with $S$.
    
    In other cases, the centralizer is the cyclic group $\mu_K$ or $\mu_{E/K}^*$, so its intersection with $S$ is cyclic and has index $1$ or $2$ in $S$, depending on whether $S$ contains elements in $\mu_K.2-\mu_K$. In these two situations, the largest normal pro-$2$ subgroup $O_2(S)$ could be either cyclic or dicyclic. The condition $g\notin \ker(\tau)$ implies $g$ is not a square in $O_2(S)$. If $O_2(S)$ is cyclic, this implies $\langle g\rangle = O_2(S)$, as in \cite{lubotzky1999lattices}. If $O_2(S)$ is dicyclic, $g$ lies in its index-$2$ core by our assumption, and we deduce that $\langle g\rangle$ equals this cyclic subgroup similarly. Since $\mathbb{F}_q$ has characteristic $2$, $C_{\slg(2,\mathbb{F}_q)}(\tau(g))\cong (C_2)^f$, implying that the centralizer $C_S(g)$ is also a $2$-group. This centralizer has order at least $\ord(g) = 2^k$, and the maximality of $O_2(S)$ confirms that $\lvert C_S(g)\rvert = 2^k$.  
\end{proof}
The next proposition concerns vertex stabilizers away from the boundary of the fixed subtree:
\begin{prop}[cf. {\cite[Proposition 4.5 (b) and (c)]{lubotzky1999lattices}}]\label{prop:not_boundary}
    Let $g\in \slg(2,K)$ be an element of order $2^k$, $k\geq 2$, with fixed tree $T_g$ and core $\mathfrak{C}_g$. Let $\Gamma$ be any lattice containing $g$, and $2^r\mid\mid \lvert \operatorname{Tor}(C_{\slg(2,K)}(g))\rvert$. Let $x$ be a vertex in $ T_g - \partial T_g$.
    
    (a) If $x$ is in $\mathfrak{C}_g$:
    \begin{itemize}
        \item If $\mathfrak{C}_g$ is a bi-infinite geodesic, one has $\Gamma_x\leq C_{(q-1)2^r}.2$.
        \item If $\mathfrak{C}_g$ is a single vertex, one has $\Gamma_x\leq C_{(q+1)2^r}.2$.
        \item If $\mathfrak{C}_g$ is a single edge, one has $\Gamma_x\leq C_{2^r}.2$, or $k = 2$ and $\Gamma_x = 2.A_4$ or $\Gamma_x = 2.S_4$.
    \end{itemize}
    
    (b) If $x$ is not in $\mathfrak{C}_g$, one has either $\Gamma_x\leq C_{2^r}.2$, or $k = 2$ and $\Gamma_x = 2.A_4$ or $\Gamma_x = 2.S_4$.
\end{prop}
\begin{proof}
    (a) Since we assume that $x\notin \partial T_g$ (which is automatic in cases (i) and (ii)), $g$ fixes the star of $x$. Realizing the full stabilizer of $x$ as $\slg(2,\mathcal{O}_K)$, we have $g\in\Gamma_{\st(x)} = \Gamma_x\cap \ker(\tau)$. Note that $\Gamma_{\st(x)}\triangleleft \Gamma_x$ and $\ker(\tau)$ is pro-$2$, so $\Gamma_x$ contains a normal $2$-subgroup containing $g$.

    By Theorem \ref{thm:finite_subgroup}, $\Gamma_x$ is contained in either a dicyclic group (cases (a) and (b)), or it is an exceptional group $2.A_4$, $2.S_4$ or $2.A_5$ (cases (c) to (e)). The first two contain a normal $2$-subgroup isomorphic to $Q_8$, while $2.A_5$ only contains $\{\pm I\}$, which is not possible here. If $\ord(g)\geq 8$, then the exceptional cases cannot hold, and an argument similar to that in \cite[Proposition 4.5]{lubotzky1999lattices} implies the assertions in cases (i) and (ii). If $\ord(g) = 4$, $\Gamma_x$ does not contain $2.A_4$ for $x\in \mathfrak{C}_g$, which is proved in Proposition \ref{prop:x_to_core}.
    
    In case (iii), while $\Gamma_x\leq C_{(q+1)2^r}.2$ (or $2.A_4$ or $2.S_4$) holds, the same argument as in \cite[Proposition 4.5]{lubotzky1999lattices} shows that, for the eigenvalues $\alpha$ and $\sigma(\alpha)$ of elements in $\Gamma_x$, where $\sigma$ generates $\operatorname{Gal}(K(\xi_{2^r})/K)$, $\alpha$ lies in the pro-$2$ group $(1+\pi_L\mathcal{O}_L) = (-1+\pi_L\mathcal{O}_L)$. It follows that $\Gamma_x\leq C_{2^r}.2$, or $\Gamma_x = 2.A_4$ or $2.S_4$ when $k=2$.
    
    (b) Take $h\in \langle g\rangle$ with $\ord(h) = 4$. Then $h$ also fixes the star $\st(x)$ and makes $O_2(\Gamma)$ larger than $\{\pm I\}$. If $\ord(g) \geq 8$, the exceptional groups $2.A_4$ and $2.S_4$ are ruled out, and the same argument as in \cite[Proposition 4.5]{lubotzky1999lattices} implies $\lvert \Gamma_x:C_{\Gamma_x}(h)\rvert = 2$. In particular, $O_2(C_{\slg(2,K)}(h))\cong C_{2^r}$, while $O_{2'}(C_{\Gamma_x}(h)) = \{1\}$ as in the proof of \cite[Proposition 4.5]{lubotzky1999lattices}. Putting these together gives the desired proposition. If $\ord(g) = 4$, the exceptional situation that $\Gamma_x = 2.A_4$ or $2.S_4$ remains possible.
\end{proof}
\begin{cor}\label{cor:ramify_no_odd}
    If $K(i)/K$ ramifies, then $\slg(2,K)$ does not contain any subgroup isomorphic to $C_{4s}$ where $s>1$ is odd.
\end{cor}
\begin{proof}
    Suppose, to the contrary, that $g,h\in C_{4s}<\slg(2,K)$, $\ord(g) = 4$ and $\ord(h) = s$. As in the proof of \cite[Proposition 4.5, (b.iii)]{lubotzky1999lattices}, writing $L = K(i)$, we have $C_G(g)\cong\{\alpha\in L^\times: N_{L/K}(\alpha)=1\}\subset (1+\pi_L\mathcal{O}_L)$ which is pro-$2$. However, $C_{4s}$ is abelian, implying $h\in C_G(g)$, a contradiction.
\end{proof}
\section{Construction of Lattices for \texorpdfstring{$\mathbb{Q}_2$}{Q2}-extensions}\label{sec:construction}
\begin{prop}[cf. {\cite[Proposition 6.1]{lubotzky1999lattices}}]
    Let $R\cong 2.A_4$, $2.S_4$ or $2.A_5$ be a finite subgroup of $\slg(2,K)$ with residue field $\mathbb{F}_q$, $q=2^f$. Suppose $R$ is contained in $\slg(2,K)_v$ for some $v\in\mathcal{V}(X)$. Let $\rho:\slg(2,K)_v\cong \slg(2,\mathcal{O}_K)\to \psl(2,\mathbb{F}_q)$ be the $\bmod \pi$ projection.
    \begin{itemize}
        \item Let $R\cong 2.A_4$. Then $\rho(R)\cong C_3$ or $A_4$. If $\rho(R)\cong C_3$, then
        \[
        \lvert \st(v)^\circ/R\rvert = \begin{cases}
            (2^f+1)/3, \quad & 2\nmid f,\\
            (2^f+5)/3, \quad & 2\mid f.
        \end{cases}
        \]
        If $\rho(R)\cong A_4$, then $2\mid f$, and
        \[
        \lvert \st(v)^\circ/R\rvert = (2^f+20)/12.
        \]
        In particular, $\lvert \st(v)^\circ/R\rvert = 1$ if and only if $q = 2$, $\lvert \st(v)^\circ/R\rvert = 2$ only if $q = 4$ and $\rho(R)\cong A_4$, and $\lvert \st(v)^\circ/R\rvert = 3$ if and only if $q=4, 8$ and $\rho(R)\cong C_3$, or $q=16$ and $\rho(R)\cong A_4$.
        \item Let $R\cong 2.S_4$. Then $\rho(R)\cong S_3$, and
        \[
        \lvert \st(v)^\circ/R\rvert = \begin{cases}
            (2^f+4)/6, \quad & 2\nmid f,\\
            (2^f+8)/6, \quad & 2\mid f.
        \end{cases}
        \]
        In particular, $\lvert \st(v)^\circ/R\rvert = 1$ if and only if $q = 2$, and $\lvert \st(v)^\circ/R\rvert = 2$ if and only if $q=4,8$.
        \item Let $R\cong 2.A_5$. Then $2\mid f$, $\rho(R)\cong A_5$, and
        \[
        \lvert \st(v)^\circ/R\rvert = \begin{cases}
            (2^f+56)/60, \quad & 4\nmid f,\\
            (2^f+104)/60, \quad & 4\mid f.
        \end{cases}
        \]
        In particular, $\lvert \st(v)^\circ/R\rvert = 1$ if and only if $q=4$, and $\lvert \st(v)^\circ/R\rvert = 2$ if and only if $q = 16,64$.
    \end{itemize}
\end{prop}
\begin{proof}
    First we determine the image $\rho(R)$. Since $\ker(\rho\mid_R)$ contains $-I$ and is normal in $R$, it is either isomorphic to $C_2$, or isomorphic to $Q_8$ when $R\cong 2.A_4$ or $2.S_4$. If $R\cong 2.A_4$, $\rho(R)\cong C_3$ or $A_4$. Since $\psl(2,\mathbb{F}_{2^f})$ does not contain an element of order $4$, it does not contain a subgroup isomorphic to $S_4$. Thus, if $R\cong 2.S_4$, we must have $\rho(R)\cong S_3$. Finally, $\rho(R)\cong A_5$ if $R\cong 2.A_5$.

    For $g\in\psl(2,\mathbb{F}_q)$, the number of leaves it fixes in $\st(v)^\circ$ is $(q+1)$ if $g = \operatorname{id}$, $1$ if $\ord(g) = 2$, $0$ if $\ord(g)\mid q+1$, and $2$ if $\ord(g)\mid q-1$. Counting the number of elements in $C_3$, $S_3$, $A_4$ and $A_5$ of each order and applying Burnside's lemma give the formulas for the numbers of orbits given above.
\end{proof}
More precisely, we can write down the exact isomorphism types of edge stabilizers:
\begin{prop}\label{prop:edge_stab}
    Suppose that $\slg(2,K)$ has residue field $\mathbb{F}_q$, $q = 2^f$, and $R<\slg(2,K)_v$ for some vertex $v$. Then for selected isomorphism types of $R$ and representatives $e$ of $\st(v)^\circ/R$, the edge stabilizer groups $G_e\cap R$ are isomorphic to the following:
    
    \begin{minipage}{0.48\textwidth}
        \begin{center}
            \begin{tabular}{|c|c|c|}
                \hline
                $R$ & $q$ & Edge Stabilizers \\
                \hline
                \multirow{4}{*}{$2.A_4$} & $2$ & $Q_8$ \\
                \cline{2-3}
                & \multirow{2}{*}{$4$} & $C_6$, $2.A_4$ \\
                \cline{3-3}
                & & $Q_8$, $2.A_4$, $2.A_4$ \\
                \cline{2-3}
                & $8$ & $Q_8$, $Q_8$, $Q_8$ \\
                \hline
                \multirow{4}{*}{$2.S_4$} & $2$ & $Q_{16}$ \\
                \cline{2-3}
                & $4$ & $Q_{16}$, $2.A_4$ \\
                \cline{2-3}
                & $8$ & $Q_8$, $Q_{16}$ \\
                \cline{2-3}
                & $16$ & $Q_{8}$, $Q_8$, $Q_{16}$, $2.A_4$ \\
                \hline
                \multirow{2}{*}{$2.A_5$} & $4$ & $2.A_4$ \\
                \cline{2-3}
                & $16$ & $C_{10}$, $2.A_4$ \\
                \hline
                $C_{4(q-1)/3}.2$ & $\geq 16$ & $C_4$, $Q_8$, $C_{4(q-1)/3}$ \\
                \hline
                $C_{2(q-1)}.2$ & $\geq 4$ & $C_4$, $C_{2(q-1)}$ \\
                \hline
                $C_{4(q-1)}$ & $\geq 4$ & $C_4$, $C_{4(q-1)}$, $C_{4(q-1)}$ \\
                \hline
                $C_{4(q-1)}.2$ & $\geq 4$ & $Q_8$, $C_{4(q-1)}$ \\
                \hline 
                $C_{4(q+1)/3}.2$ & $\geq 8$ & $C_4$, $Q_8$ \\
                \hline
                $C_{2(q+1)}.2$ & & $C_4$ \\
                \hline
                $C_{4(q+1)}$ & & $C_4$ \\
                \hline
                $C_{4(q+1)}.2$ & & $Q_8$ \\
                \hline
            \end{tabular}
        \end{center}
    \end{minipage}
    \hfill
    \begin{minipage}{0.48\textwidth}
        \begin{center}
            \begin{tabular}{|c|c|c|}
                \hline
                $R$ & $q$ & Edge Stabilizers \\
                \hline
                \multirow{7}{*}{$Q_{2^k}$} & \multirow{3}{*}{$2$} & $Q_{2^{k-1}}$, $Q_{2^k}$ \\
                \cline{3-3}
                & & $C_{2^{k-1}}$, $Q_{2^k}$ \\
                \cline{3-3}
                & & $Q_{2^k}$, $Q_{2^k}$, $Q_{2^k}$\\
                \cline{2-3}
                & \multirow{4}{*}{$4$} & $C_{2^{k-2}}$, $Q_{2^k}$ \\
                \cline{3-3}
                & & $C_{2^{k-1}}$, $C_{2^{k-1}}$, $Q_{2^k}$ \\
                \cline{3-3}
                & & $Q_{2^{k-1}}$, $Q_{2^{k-1}}$, $Q_{2^k}$ \\
                \cline{3-3}
                & & $Q_{2^k}$, $Q_{2^k}$, $Q_{2^k}$, $Q_{2^k}$, $Q_{2^k}$\\
                \hline
                \multirow{4}{*}{$C_{2^k}$} & \multirow{2}{*}{$2$} & $C_{2^{k-1}}$, $C_{2^k}$ \\
                \cline{3-3}
                & & $C_{2^k}$, $C_{2^k}$, $C_{2^k}$\\
                \cline{2-3}
                & \multirow{2}{*}{$4$} & $C_{2^{k-1}}$, $C_{2^{k-1}}$, $C_{2^k}$ \\
                \cline{3-3}
                & & $C_{2^k}$, $C_{2^k}$, $C_{2^k}$, $C_{2^k}$, $C_{2^k}$\\
                \hline
                \multirow{2}{*}{$C_{3\cdot 2^k}.2$} & $2$ & $Q_{2^{k+1}}$ \\
                \cline{2-3}
                & $4$ & $Q_{2^{k+1}}$, $C_{3\cdot 2^k}$ \\
                \hline
                \multirow{2}{*}{$C_{3\cdot 2^k}$} & $2$ & $C_{2^{k}}$ \\
                \cline{2-3}
                & $4$ & $C_{2^{k}}$, $C_{3\cdot 2^k}$, $C_{3\cdot 2^k}$ \\
                \hline
                $C_{5\cdot 2^k}.2$ & $4$ & $Q_{2^{k+1}}$ \\
                \hline
                $C_{5\cdot 2^k}$ & $4$ & $C_{2^{k}}$ \\
                \hline
            \end{tabular}
        \end{center}
    \end{minipage}
    Here, $C_{4(q-1)/3}.2$ requires $2\mid f$, and $C_{4(q+1)/3}.2$ requires $2\nmid f$.
\end{prop}
The proof of Proposition \ref{prop:edge_stab} is left to the reader.

\textbf{Examples of Lattices in $\slg(2,K)$ for finite extensions $K/\mathbb{Q}_2$.}

\begin{exm}\label{exm:2}
    Suppose $K = \mathbb{Q}_2(\xi_{q-1})$, so the residue field is $\mathbb{F}_q$ and $r_2 = 1$. We claim that there exists a lattice $\Gamma<\slg(2,K)$, with the following graph of groups
    \begin{center}
        \begin{tikzpicture}[scale=2.5,every node/.style={font=\footnotesize}]
            \draw[thick] (0,0) --node[above]{$C_4$} (1,0);
            \fill[] (0,0) circle(0.025) node[left]{$C_{2(q+1)}.2$};
            \fill[] (1,0) circle(0.025) node[right]{$C_{2(q+1)}.2$};
        \end{tikzpicture}
    \end{center}
    Indeed, fix a vertex $u\in \mathcal{V}(T)$, and a subgroup $\Gamma_u<G_u$ such that $\Gamma_u\cong C_{2(q+1)}.2$. Its subgroup of order $q+1$ acts transitively on $\st(u)$, so for any edge $e$ incident with $u$, we have $\Gamma_e\coloneqq \Gamma_u\cap G_e\cong C_4$. Since $K(i)/K$ is ramified and $r_2=1$, the fixed tree of every subgroup isomorphic to $C_4$ is a single edge. In particular, $T_{\Gamma_e} = e$, and its other endpoint satisfies $v\in \partial T_{\Gamma_e}$. 
    
    By Proposition \ref{prop:cyclic_conj}, there is a unique $H_v$-conjugacy class of subgroups of $G_v$ isomorphic to $C_4$. For a chosen subgroup of $G_v$ isomorphic to $C_{2(q+1)}.2$, \cite[Lemma 4.7(2)]{lubotzky1999lattices} guarantees an $H_v$-conjugate $\Gamma_v<G_v$ with $\Gamma_v>\Gamma_e$. Since $\Gamma_v\cong C_{2(q+1)}.2$, $\Gamma_v\cap G_e\cong C_4$, so $\Gamma_v\cap G_e = \Gamma_e$ as well. \cite[Corollary 5.2 and Theorem 5.4]{lubotzky1999lattices} then imply the existence of the desired lattice $\Gamma$. The covolume of this lattice is $\frac{1}{4(q+1)}$.
\end{exm}
\begin{exm}[cf. {\cite[Example (6.3)]{lubotzky1999lattices}}]\label{exm:3}
    Let $K = \mathbb{Q}_2(\xi_3)$ ($q=4$, $e(K/\mathbb{Q}_2) = 1$). We claim that there exists a lattice $\Gamma<\slg(2,K)$, with the following graph of groups
    \begin{center}
        \begin{tikzpicture}[scale=2.5,every node/.style={font=\footnotesize}]
            \draw[thick] (0,0) --node[above]{$2.A_4$} (1,0);
            \fill[] (0,0) circle(0.025) node[left]{$2.A_5$};
            \fill[] (1,0) circle(0.025) node[right]{$2.A_5$};
        \end{tikzpicture}
    \end{center}
    Indeed, fix a vertex $u\in \mathcal{V}(T)$ and a subgroup $\Gamma_u<G_u$ such that $\Gamma_u\cong 2.A_5$. Since $2.A_5$ acts transitively on $\st(u)$, for any edge $e$ incident with $u$, we have $\Gamma_e\coloneqq \Gamma_u\cap G_e\cong 2.A_4$. Since $r_2 = 1$, Corollary \ref{cor:A4_attach} implies that $T_{\Gamma_e}$ is a path of length $1$, so the other endpoint satisfies $v\in\partial T_{\Gamma_e}$. As in the previous example, Corollary \ref{cor:A4_attach} implies the existence of a subgroup $\Gamma_v<G_v$ with $\Gamma_v\cong 2.A_5$ and $\Gamma_v\cap G_e = \Gamma_e$. \cite[Corollary 5.2 and Theorem 5.4]{lubotzky1999lattices} then imply the existence of the desired lattice $\Gamma$. The covolume of this lattice is $\frac{1}{120}$.
\end{exm}
\begin{exm}[cf. {\cite[Example (6.4)]{lubotzky1999lattices}}]\label{exm:4}
    Let $K = \mathbb{Q}_2(\xi_{15})$ ($q=16$, $e(K/\mathbb{Q}_2) = 1$). We claim that there exists a lattice $\Gamma<\slg(2,K)$, with the following graph of groups
    \begin{center}
        \begin{tikzpicture}[scale=2.5,every node/.style={font=\footnotesize}]
            \draw[thick] (0,0) to[out=45, in=135]node[above]{$2.A_4$} (1,0);
            \draw[thick] (0,0) to[out=-45, in=-135]node[below]{$C_{10}$} (1,0);
            \fill[] (0,0) circle(0.025) node[left]{$2.A_5$};
            \fill[] (1,0) circle(0.025) node[right]{$2.A_5$};
        \end{tikzpicture}
    \end{center}
    As before, choose $u\in\mathcal{V}(T)$ and $\Gamma_u<G_u$ such that $\Gamma_u\cong 2.A_5$. Let $e_1 = [u,v]$ be an edge in $\st(u)$ such that $\Gamma_{e_1}\coloneqq \Gamma_u\cap G_{e_1}\cong 2.A_4$. Similarly, $r_2 = 1$ implies there is $\Gamma_v<G_v$ with $\Gamma_v\cong 2.A_5$ and $\Gamma_v\cap G_{e_1} = \Gamma_{e_1}$. Let $e_2 = [v,w]$ be an edge in the other $\Gamma_v$-orbit of $\st(v)$. Then $\Gamma_{e_2}\coloneqq \Gamma_v\cap G_{e_2}\cong C_{10}$.

    Now, since $d(u,w) = 2$, there is $g\in G$ such that $g.u = w$. By Corollary \ref{cor:cyclic}, there is a unique $G_w$-conjugacy class of $G_w$-subgroups isomorphic to $C_{10}$. Then \cite[Lemma 4.7 (1)]{lubotzky1999lattices} implies there is $\Gamma_w<G_w$ among the conjugates $g\Gamma_u g^{-1}<G_w$ for such $g$, such that $\Gamma_w\cap G_{e_2} = \Gamma_{e_2}$. Using this element $g$ to define the boundary monomorphism, we construct the claimed graph of groups. \cite[Corollary 5.2 and Theorem 5.4]{lubotzky1999lattices} then imply the existence of the desired lattice $\Gamma$. The covolume of this lattice is $\frac{1}{120}$.
\end{exm}
\begin{exm}[cf. {\cite[Example (6.5)]{lubotzky1999lattices}}]\label{exm:a}
    Let $K = \mathbb{Q}_2(z_8) = \mathbb{Q}_2(\sqrt{2})$ ($q=2$, $e(K/\mathbb{Q}_2) = 2$). In this case, $\slg(2,K)$ contains elements of order $8$, and the fixed tree of any subgroup isomorphic to $Q_{16}$ has diameter $r_3 = r_2/2 = 1$, i.e., it is an edge. Using Proposition \ref{prop:supergroup_exceptional} and arguing as in Example \ref{exm:4}, we obtain a lattice in $\slg(2,K)$ with the following graph of groups
    \begin{center}
        \begin{tikzpicture}[scale=2.5,every node/.style={font=\footnotesize}]
            \draw[thick] (0,0) --node[above]{$Q_{16}$} (1,0);
            \fill[] (0,0) circle(0.025) node[left]{$2.S_4$};
            \fill[] (1,0) circle(0.025) node[right]{$2.S_4$};
        \end{tikzpicture}
    \end{center}
    The covolume of this lattice is $\frac{1}{48}$.
\end{exm}
\begin{exm}[cf. {\cite[Example (6.8)]{lubotzky1999lattices}}]\label{exm:5}
    Let $K = \mathbb{Q}_2(z_{12}) = \mathbb{Q}_2(\sqrt{3})$ ($q=2$, $e(K/\mathbb{Q}_2) = 2$). We claim that there exists a lattice $\Gamma<\slg(2,K)$, with the following graph of groups
    \begin{center}
        \begin{tikzpicture}[scale=2.5,every node/.style={font=\footnotesize}]
            \draw[thick] (0,0) --node[above]{$Q_8$} (1,0);
            \fill[] (0,0) circle(0.025) node[left]{$2.A_4$};
            \fill[] (1,0) circle(0.025) node[right]{$C_{12}.2$};
        \end{tikzpicture}
    \end{center}
    As above, as $\lvert K/\mathbb{Q}\rvert = 2$, we can take a vertex $u$ and $\Gamma_u<G_u$, $\Gamma_u\cong 2.A_4$; for an edge $e=[u,v]$, $\Gamma_e = \Gamma_u\cap G_e\cong Q_8$. $\Gamma_u\cong 2.A_4$ implies $u$ is the center of $T_{\Gamma_e}$, and since $r_2 = 2$, $v\in\partial T_{\Gamma_e}$. Since $K(i)/K$ is unramified, Proposition \ref{prop:dicyclic_subgroups} (2) implies there is a subgroup of $G_u$ isomorphic to $C_{12}.2$. Since the three cores of order-$4$ elements in $\Gamma_e$ are distinct and at distance $r_2/2=1$ from $u$, these exhaust the vertices in $\partial T_{\Gamma_e}$. Proposition \ref{prop:supergroup} implies the existence of $\Gamma_v<G_v$, $\Gamma_v\cong C_{12}.2$ and $\Gamma_v\cap G_e = \Gamma_e$. \cite[Corollary 5.2 and Theorem 5.4]{lubotzky1999lattices} then imply the existence of the desired lattice $\Gamma$. The covolume of this lattice is $\frac{1}{24}$.
\end{exm}
\begin{exm}\label{exm:6}
    Let $q=2$, $e(K/\mathbb{Q}_2) = 2$, and $K\neq \mathbb{Q}_2(z_{12})$ or $\mathbb{Q}_2(\xi_4)$, so $K(i)/K$ is ramified. We claim that there exists a lattice $\Gamma<\slg(2,K)$, with the following graph of groups
    \begin{center}
        \begin{tikzpicture}[scale=2.5,every node/.style={font=\footnotesize}]
            \draw[thick] (0,0) --node[above]{$Q_8$} (1,0) --node[above]{$C_4$} (2,0);
            \fill[] (0,0) circle(0.025) node[below]{$2.A_4$};
            \fill[] (1,0) circle(0.025) node[below]{$Q_8$};
            \fill[] (2,0) circle(0.025) node[below]{$C_{6}.2$};
        \end{tikzpicture}
    \end{center}
    Again, we take $\Gamma_u\cong 2.A_4$, $e_1 = [u,v]$, $\Gamma_{e_1}\cong Q_8$; in this case, we take $\Gamma_v = \Gamma_{e_1}$. There are two $\Gamma_v$-orbits of $\st(v)$; the two edges other than $e_1$ are in the same orbit. Let $e_2 = [v,w]$ be one of them. Then $\Gamma_{e_2}\cong C_4$. Let $g$ be a generator of $\Gamma_{e_2}$. Since $K(i)/K$ is ramified and $r_2 = 2$, Theorem \ref{thm:width_formula} implies $\mathfrak{C}_g$ is an edge and $T_g = N(\mathfrak{C}_g,1)$. Since $\mathfrak{C}_g\subset T_{\Gamma_{e_1}}$ and $w\in \mathcal{V}(T_g)$, $\mathfrak{C}_g = e_1$ and $w\in\partial T_g$. \cite[Proposition 4.7(2)]{lubotzky1999lattices} and Proposition \ref{prop:cyclic_conj} then imply the existence of $\Gamma_w<G_w$, $\Gamma_w\cong C_6.2$ and $\Gamma_w\cap G_{e_2} = \Gamma_{e_2}$. \cite[Corollary 5.2 and Theorem 5.4]{lubotzky1999lattices} then imply the existence of the desired lattice $\Gamma$. The covolume of this lattice is $\frac{1}{8}$.
\end{exm}
\begin{exm}[cf. \cite{lubotzky1999lattices}, Example 6.9]\label{exm:x}
    Let $K = \mathbb{Q}_2(\xi_4)$ ($q=2$, $e(K/\mathbb{Q}_2)=2$). Similarly to Example \ref{exm:6}, we start with an edge with group $Q_8$, such that the two endpoints have groups $2.A_4$ and $Q_8$, respectively, and an additional $C_4$-edge is incident with the $Q_8$ vertex. In this case, since $i = \xi_4\in K$, Proposition \ref{prop:dicyclic_fixed} implies that the core of the generator of this $C_4$ is a bi-infinite path meet the $Q_8$-fixed tree at a vertex. A hyperbolic element $g\in \gl(2,K)$ along this axis with translation length $1$ necessarily centralize this $C_4$-subgroup; conjugate the $Q_8$-edge configuration by $g$, we obtain a lattice with the following graph of groups:
    \begin{center}
        \begin{tikzpicture}[scale=2.5,every node/.style={font=\footnotesize}]
            \draw[thick] (0,0) --node[above]{$Q_8$} (1,0) --node[above]{$C_4$} (2,0)--node[above]{$Q_8$} (3,0);
            \fill[] (0,0) circle(0.025) node[below]{$2.A_4$};
            \fill[] (1,0) circle(0.025) node[below]{$Q_8$};
            \fill[] (2,0) circle(0.025) node[below]{$Q_8$};
            \fill[] (3,0) circle(0.025) node[below]{$2.A_4$};
        \end{tikzpicture}
    \end{center}
    The covolume of this lattice is $\frac{1}{6}$.
\end{exm}
\begin{exm}\label{exm:7}
    Let $K = \mathbb{Q}_2(z_{20})$ ($q=4$, $e(K/\mathbb{Q}_2) = 2$). We claim that there exists a lattice $\Gamma<\slg(2,K)$, with the following graph of groups
    \begin{center}
        \begin{tikzpicture}[scale=2.5,every node/.style={font=\footnotesize}]
            \draw[thick] (0,0) --node[above]{$2.A_4$} (1,0) --node[above]{$2.A_4$} (2,0);
            \draw[thick] (1,0) --node[left]{$Q_8$} (1,1);
            \fill[] (0,0) circle(0.025) node[below]{$2.A_5$};
            \fill[] (1,0) circle(0.025) node[below]{$2.A_4$};
            \fill[] (2,0) circle(0.025) node[below]{$2.A_5$};
            \fill[] (1,1) circle(0.025) node[above]{$C_{20}.2$};
        \end{tikzpicture}
    \end{center}
    Similarly, take a subgroup $\Gamma_u\cong 2.A_4$, such that $u$ is the center of $T_{\Gamma_u}$; in this case, $\st(u)$ has three $\Gamma_u$-orbits; two of them contain a single edge, denoted by $e_i = [u,v_i]$, $i=1,2$, so $\Gamma_{e_i} = \Gamma_u$. Since $r_2 = 2$, $v_1$, $v_2\in\partial T_{\Gamma_u}$, and we can choose $\Gamma_{v_i}\cong 2.A_5$. An edge $e_3 = [u,w]$ in the other orbit has edge group $\Gamma_{e_3} \cong Q_8$; Lemma \ref{lem:exceptional_tree_structure} implies $u$ is the center of $T_{\Gamma_{e_3}}$. Since $K(i)/K$ is unramified, the three vertices in this orbit are the cores of the three order-$4$ subgroups of $\Gamma_{e_3}$. Proposition \ref{prop:supergroup} therefore allows us to choose $\Gamma_w\cong C_{20}.2$ with $\Gamma_w\cap G_{e_3}=\Gamma_{e_3}$. This completes the construction of the prescribed graph of groups, whose covolume is $\frac{1}{24}$.
\end{exm}
\begin{exm}\label{exm:b}
    Let $K = \mathbb{Q}_2(z_5,z_8)$ ($q=4$, $e(K/\mathbb{Q}_2) = 2$). In this case, $\slg(2,K)$ contains both $2.S_4$ and $2.A_5$ subgroups. Since $r_3 = 1$, we can first construct two $2.S_4$ vertices connected by a $Q_{16}$ edge using Proposition \ref{prop:supergroup_exceptional}, as in Example \ref{exm:a}. The other edge groups are $2.A_4$, and by Corollary \ref{cor:A4_attach}, the other endpoints of these edges can have vertex groups isomorphic to $2.A_5$. We obtain a lattice in $\slg(2,K)$ with the following graph of groups:
    \begin{center}
        \begin{tikzpicture}[scale=2.5,every node/.style={font=\footnotesize}]
            \draw[thick] (0,0) --node[above]{$2.A_4$} (1,0) --node[above]{$Q_{16}$} (2,0) --node[above]{$2.A_4$} (3,0);
            \fill[] (0,0) circle(0.025) node[below]{$2.A_5$};
            \fill[] (1,0) circle(0.025) node[below]{$2.S_4$};
            \fill[] (2,0) circle(0.025) node[below]{$2.S_4$};
            \fill[] (3,0) circle(0.025) node[below]{$2.A_5$};
        \end{tikzpicture}
    \end{center}
    The covolume of this lattice is $\frac{7}{240}$.
\end{exm}
\begin{exm}\label{exm:c}
    Let $K = \mathbb{Q}_2(z_8,z_{12})$ ($q=2$, $e(K/\mathbb{Q}_2) = 4$). In this case, $\slg(2,K)$ contains a $Q_{16}$ subgroup. Since $r_3=2$, if we begin with the central vertex of $T_{Q_{16}}$, there are three edge orbits of its star, and all of them have edge groups $Q_{16}$. Proposition \ref{prop:supergroup} implies that the core of the order-$8$ element in $Q_{16}$ allows a $C_{24}.2$ group, while Proposition \ref{prop:supergroup_exceptional} implies that the centers of the fixed subtrees of the two $Q_8$ subgroups allow $2.S_4$ groups; all these are distinct vertices on $\partial T_{Q_{16}}$. We obtain a lattice in $\slg(2,K)$ with the following graph of groups:
    \begin{center}
        \begin{tikzpicture}[scale=2.5,every node/.style={font=\footnotesize}]
            \draw[thick] (0,0) --node[above]{$Q_{16}$} (1,0) --node[above]{$Q_{16}$} (2,0);
            \draw[thick] (1,0) --node[left]{$Q_{16}$} (1,1);
            \fill[] (0,0) circle(0.025) node[below]{$2.S_4$};
            \fill[] (1,0) circle(0.025) node[below]{$Q_{16}$};
            \fill[] (2,0) circle(0.025) node[below]{$2.S_4$};
            \fill[] (1,1) circle(0.025) node[above]{$C_{24}.2$};
        \end{tikzpicture}
    \end{center}
    The covolume of this lattice is $\frac{1}{16}$.
\end{exm}
\begin{exm}\label{exm:d}
    Let $K = \mathbb{Q}_2(z_{16})$ ($q=2$, $e(K/\mathbb{Q}_2) = 4$). In this case, $\slg(2,K)$ contains a $Q_{32}$ subgroup that fixes only an edge $[u,v]$, and Theorem \ref{thm:width_formula} implies that the $C_{16}$ subgroup fixes exactly the same edge. Therefore, a representative of the other $Q_{32}$-orbit of $\st(u)^\circ$ has edge group $Q_{16}$, and Proposition \ref{prop:supergroup_exceptional} guarantees a vertex group isomorphic to $2.S_4$. Applying the same argument at $v$, we obtain a lattice in $\slg(2,K)$ with the following graph of groups:
    \begin{center}
        \begin{tikzpicture}[scale=2.5,every node/.style={font=\footnotesize}]
            \draw[thick] (0,0) --node[above]{$Q_{16}$} (1,0) --node[above]{$Q_{32}$} (2,0) --node[above]{$Q_{16}$} (3,0);
            \fill[] (0,0) circle(0.025) node[below]{$2.S_4$};
            \fill[] (1,0) circle(0.025) node[below]{$Q_{32}$};
            \fill[] (2,0) circle(0.025) node[below]{$Q_{32}$};
            \fill[] (3,0) circle(0.025) node[below]{$2.S_4$};
        \end{tikzpicture}
    \end{center}
    The covolume of this lattice is $\frac{5}{96}$.
\end{exm}
\begin{exm}\label{exm:y}
    Let $K = \mathbb{Q}_2(\xi_{12})$ ($q=4$, $e(K/\mathbb{Q}_2) = 2$). Similarly to Example \ref{exm:7}, we begin with a vertex with group $2.A_4$, with two $2.A_4$ edges incident with $2.A_5$ vertices and an additional $Q_8$ edge. In this case, a $C_{12}.2$ vertex is allowed instead as the other endpoint of the $Q_8$ edge. Since $i\in K$ and $3\mid (q-1)$, the order-$12$ element in this dicyclic group fixes a bi-infinite path, attaching to the $Q_8$-tree at the current vertex, providing an outgoing $C_{12}$-edge. A conjugation of the entire configuration similar to Example \ref{exm:x} gives us a lattice with the following graph of groups:
    \begin{center}
        \begin{tikzpicture}[scale=2.5,every node/.style={font=\footnotesize}]
            \draw[thick] (0,0) --node[above]{$2.A_4$} (1,0) --node[above]{$Q_8$} (2,0) --node[above]{$C_{12}$} (3,0)--node[above]{$Q_8$} (4,0)--node[above]{$2.A_4$} (5,0);
            \draw[thick] (1,0) --node[left]{$2.A_4$}(1,1);
            \draw[thick] (4,0) --node[left]{$2.A_4$}(4,1);
            \fill[] (0,0) circle(0.025) node[below]{$2.A_5$};
            \fill[] (1,0) circle(0.025) node[below]{$2.A_4$};
            \fill[] (2,0) circle(0.025) node[below]{$C_{12}.2$};
            \fill[] (3,0) circle(0.025) node[below]{$C_{12}.2$};
            \fill[] (4,0) circle(0.025) node[below]{$2.A_4$};
            \fill[] (5,0) circle(0.025) node[below]{$2.A_5$};
            \fill[] (1,1) circle(0.025) node[above]{$2.A_5$};
            \fill[] (4,1) circle(0.025) node[above]{$2.A_5$};
        \end{tikzpicture}
    \end{center}
    The covolume of this lattice is $\frac{1}{10}$.
\end{exm}
\begin{exm}\label{exm:z}
    Let $K = \mathbb{Q}_2(\xi_{8})$ ($q=2$, $e(K/\mathbb{Q}_2) = 4$). Similarly to Example \ref{exm:c}, we begin with a $Q_{16}$-tree with radius $1$, with two of the endpoints assigned to $2.S_4$ groups. In this case, we attach a $Q_{16}$-vertex to the remaining edge. Since $\xi_8\in K$, by Proposition \ref{prop:dicyclic_fixed}, the core of the order-$8$ element in this $Q_{16}$ is a bi-infinite path, and by Lemma \ref{lem:exceptional_tree_structure}, this core is attached to the vertex different from the two $2.S_4$-vertices. Therefore the $C_8$-axis gives an outgoing $C_8$-edge from the $Q_{16}$-vertex above, and a conjugation of the configuration similar to Example \ref{exm:y} gives birth to a lattice with the following graph of groups: 
    \begin{center}
        \begin{tikzpicture}[scale=2.5,every node/.style={font=\footnotesize}]
            \draw[thick] (0,0) --node[above]{$Q_{16}$} (1,0) --node[above]{$Q_{16}$} (2,0) --node[above]{$C_8$} (3,0)--node[above]{$Q_{16}$} (4,0)--node[above]{$Q_{16}$} (5,0);
            \draw[thick] (1,0) --node[left]{$Q_{16}$}(1,1);
            \draw[thick] (4,0) --node[left]{$Q_{16}$}(4,1);
            \fill[] (0,0) circle(0.025) node[below]{$2.S_4$};
            \fill[] (1,0) circle(0.025) node[below]{$Q_{16}$};
            \fill[] (2,0) circle(0.025) node[below]{$Q_{16}$};
            \fill[] (3,0) circle(0.025) node[below]{$Q_{16}$};
            \fill[] (4,0) circle(0.025) node[below]{$Q_{16}$};
            \fill[] (5,0) circle(0.025) node[below]{$2.S_4$};
            \fill[] (1,1) circle(0.025) node[above]{$2.S_4$};
            \fill[] (4,1) circle(0.025) node[above]{$2.S_4$};
        \end{tikzpicture}
    \end{center}
    The covolume of this lattice is $\frac{1}{6}$.
\end{exm}
\section{Proof of Minimal Covolume Results}\label{sec:proof}
\begin{prop}[cf. {\cite[Proposition C]{lubotzky1999lattices}}]\label{prop:no_8}
    For a finite extension $K/\mathbb{Q}_2$ with residue field of order $q=2^f$, denote by $\lambda'(K)$ the minimal covolume among cocompact lattices $\Gamma<\slg(2,K)$ without elements of order $8$. Then:
    \[
    \lambda'(K) = \begin{cases}
        \frac{1}{120}, & \ K = \mathbb{Q}_2(\xi_3) \text{ or }\mathbb{Q}_2(\xi_{15}), \\
        \frac{1}{24}, & \ K = \mathbb{Q}_2(z_{12}) = \mathbb{Q}_2(\sqrt{3}),\text{ or }K = \mathbb{Q}_2(\xi_3,z_{20}),\\
        \frac{1}{8}, & \ q=2,\ e(K/\mathbb{Q}_2) = 2,\text{ and }K\neq \mathbb{Q}_2(z_{12}),\ \mathbb{Q}_2(\xi_4),\\
        \frac{1}{10}, & \ K = \mathbb{Q}_2(\xi_{12}), \\
        \frac{1}{6}, & \ K = \mathbb{Q}_2(\xi_4), \\
        \frac{1}{4(q+1)}, & \ K = \mathbb{Q}_2(\xi_{q-1}),\ q\ne 4,\ 16, \\
        \frac{1}{2(q+1)}, & \ \text{in all other cases.}
    \end{cases}
    \]
\end{prop}
\begin{proof}
    The upper bound $\lambda'(K)\le \frac{1}{2(q+1)}$ follows directly from Example \ref{exm:reg} and \cite[proof of Proposition C]{lubotzky1999lattices}. For each $q = 2^f$, we discuss below only lattices with smaller or equal covolumes.
    
    \textbf{Case (1): $q>2^4$}. The binary polyhedral groups $2.A_4$ and $2.A_5$ cannot occur in the minimal covolume lattices for $q\geq 2^5$: when $q=2^5$, only $2.A_4$ is admitted and $\frac{1}{24}>\frac{1}{2(2^5+1)}$, and when $q\geq 2^6$, $\frac{1}{120}>\frac{1}{2(2^6+1)}$. Since there are no elements of order $8$, Theorem \ref{thm:finite_subgroup} implies any vertex group is a subgroup of $C_{4(q+1)}.2$ or $C_{4(q-1)}.2$. If a vertex group of $C_{2(q+1)}\cong C_{(q+1)}.2$ exists, then the unique $C_2$-edge implies that Example \ref{exm:reg} is the only lattice with covolume $\frac{1}{2(q+1)}$. For other graphs of groups, given that $\lambda'(K)\le \frac{1}{2(q+1)}$, only the following vertex groups are allowed:
    \begin{itemize}
        \item $C_{4(q+1)/3}.2$, $C_{2(q+1)}.2$, $C_{4(q+1)}$, and $C_{4(q+1)}.2$, and
        \item $C_{4(q-1)/3}.2$, $C_{2(q-1)}.2$, $C_{4(q-1)}$, and $C_{4(q-1)}.2$.
    \end{itemize}
    Here, if a vertex group is isomorphic to $C_{4(q\pm1)/3}.2$ or $C_{4(q\pm1)}.2$, then the stabilizers of some edges in its star are isomorphic to $Q_8$, and the group at the other endpoint is also of one of these isomorphism types. In these cases, Proposition \ref{prop:not_boundary} implies both endpoints are in $\partial T_{Q_8}$, so Proposition \ref{prop:dicyclic_fixed} implies that $r_2 = 1$ and $K(i)/K$ is ramified. However, as $q>4$, Corollary \ref{cor:ramify_no_odd} implies $\slg(2,K)$ does not contain these subgroups, a contradiction.

    Now, $C_{2(q\pm 1)}.2$ and $C_{4(q\pm 1)}$ become the only vertex group candidates. Since $\lambda'(K)\le \frac{1}{2(q+1)}$, a vertex group isomorphic to $C_{2(q - 1)}.2$ or $C_{4(q - 1)}$ prevents another vertex group of the same color from existing; since the remaining candidate vertex groups have order at most $4(q+1)$ and $\frac{2}{4(q+1)} > \frac{1}{4(q-1)}$, there is also only one vertex of the opposite color (of the same order $4(q-1)$), making the graph have only two vertices. For a graph made from only vertex groups isomorphic to $C_{2(q+1)}.2$ or $C_{4(q+1)}$, since they act transitively on the star, the graph of groups can have only two vertices as well. If both vertex groups are cyclic groups isomorphic to $C_{4(q\pm 1)}$, the subgroup in $\slg(2,K)$ is uniquely determined by the edge group isomorphic to $C_4$; the two vertex groups are identical, a contradiction. The $C_{4(q+1)}*_{C_4} C_{2(q+1)}.2$ case is not possible either: suppose $g\in C_4$ is the generator. Since $C_{4(q+1)}$ is abelian, the corresponding vertex is not in $\partial T_g$, so Proposition \ref{prop:not_boundary} implies $K(i)/K$ is quadratic and unramified, and this vertex is $\mathfrak{C}_g$. Then the vertex for $C_{2(q+1)}.2$ should be in $\partial T_g$ and $r_2 = 1$, a contradiction. 
    
    This leaves only two possibilities, namely the graphs of groups with both vertex groups isomorphic to $C_{2(q\pm 1)}.2$. For the $C_4$ edge in both cases, this $C_4$ subgroup in $C_{2(q\pm 1)}.2$ acts on the star with order $2$, so it fixes only one edge in the star. Thus $T_{C_4}$ has only one edge, and $r_2 = 1$. Therefore, when $r_2 = 1$, Example \ref{exm:2} gives $\lambda'(K) = \frac{1}{4(q+1)}$, while $\lambda'(K) = \frac{1}{2(q+1)}$ holds when $r_2>1$. 

    \textbf{Case (2): $q = 2^4$.} In this case, the vertex group $2.A_5$ is allowed, and an edge stabilizer is $2.A_4$; the other vertex group has to be $2.A_5$ as well. Since $2.A_4$ is self-normalizing in $2.A_5$ and the other edge stabilizer is $C_{10}$, this $2.A_4$ fixes only this connecting edge. Hence the fixed subtree $T_{2.A_4}$ contains only one edge, which forces $r_2 = 1$. Example \ref{exm:4} gives $\lambda'(K) = \frac{1}{120}$ when $r_2 = 1$, and combining this with the argument in Case (1), we see that $\lambda'(K) = \frac{1}{2(q+1)}$ holds when $r_2>1$.

    \textbf{Case (3): $q = 2^3$.} In this case, the vertex group $\Gamma_x\cong 2.A_4$ is allowed, and the corresponding star has three $2.A_4$-orbits, all with edge groups isomorphic to $Q_8$; it follows that these $Q_8$ subgroups all coincide with the unique such subgroup in $2.A_4$. Since $Q_8$ determines a unique supergroup isomorphic to $2.A_4$, the vertex stabilizer at the other endpoint is one of $C_{12}.2$, $C_{28}.2$, or $C_{36}.2$. Therefore, all the nine boundary vertices of $\st(x)$ are fixed by this $Q_8$ subgroup and have vertex stabilizers $C_{4s}.2$ with odd $s>1$. However, Proposition \ref{prop:supergroup} allows at most three vertex stabilizers of this form, a contradiction. Therefore the conclusion remains the same as in Case (1).

    \textbf{Case (4): $q = 2^2$.} In this case, $C_{4(q-1)/3}.2 = Q_8$ no longer satisfies the minimal covolume requirement as a vertex group; in particular, the same argument as in Case (1) still holds. Now if the exceptional groups $2.A_4$ and $2.A_5$ occur, there is a $2.A_4$ subgroup, and Corollary \ref{cor:A4_attach} implies its fixed subtree is a path of length $r_2$. When $r_2 = 1$, Example \ref{exm:3} gives the unique graph of groups for lattices of minimal covolume $\frac{1}{120}$. 
    
    When $r_2>1$, we claim that the fixed subtree of $2.A_4$ appears in its entirety: indeed, without elements of order $8$, $N_\Gamma(2.A_4)$ is $2.A_4$ itself; for two vertices $x,y$ in the same $\Gamma$-orbit of $T_{2.A_4}$, identified by $\gamma\in \Gamma$, the conjugate of $2.A_4$ by $\gamma$ is identified with $2.A_4$, up to an additional conjugation in $2.A_5$. Therefore $\gamma$ fixes the vertex $x=y$. Now if $r_2 \ge 4$, then the center of $T_{2.A_4}$ is attached to two distinct vertex groups isomorphic to $2.A_4$; the third edge orbit has edge group $Q_8$ and again cannot attach to another $2.A_4$ group, so it requires an additional vertex group isomorphic to $C_{12}.2$ or $C_{20}.2$. The resulting covolume is
    \[
    \mu(\Gamma\backslash G)\geq \frac{1}{24}+ \frac{1}{24} + \frac{1}{40}>\frac{1}{10},
    \]
    so it is not of minimal covolume.

    Now consider the cases $r_2 = 2$ or $3$. If $r_2 = 3$, then $K(i)/K$ is ramified; again as the $Q_8$-edge group requires $C_{12}.2$ or $C_{20}.2$, this contradicts Corollary \ref{cor:ramify_no_odd}. If $r_2 = 2$ and $\slg(2,K)$ contains a subgroup isomorphic to $C_{20}.2$, then $K = \mathbb{Q}_2(z_{20})$ and Example \ref{exm:7} gives the minimal covolume $\lambda'(K) = \frac{1}{24}$. Otherwise, the vertex group incident with the $Q_8$-edge group has to be $C_{12}.2$. If the other endpoint of the $C_{12}$ edge is $C_{12}$, then the covolume is at least $\frac{1}{24}+\frac{1}{12}>\frac{1}{10}$.
    
    Therefore, the $C_{12}$ edge is incident with another $C_{12}.2$ group, leading to a path in the graph of groups with endpoint groups $2.A_4$ and intervening vertex groups $C_{12}.2$; there is an even number of such intervening vertices as the edge stabilizers alternate between $C_{12}$ and $Q_8$. Below is the smallest possible configuration with two $C_{12}.2$ vertices:
    \begin{center}
        \begin{tikzpicture}[scale=2.5,every node/.style={font=\footnotesize}]
            \draw[thick,dashed] (-0.5,0)node[above]{$2.A_4$} -- (0,0) --node[above]{$Q_8$} (1,0) --node[above]{$C_{12}$} (2,0) --node[above]{$Q_8$} (3,0) -- (3.5,0)node[above]{$2.A_4$};
            \draw[thick,dashed] (0,0) -- (0,0.5)node[left]{$2.A_4$};
            \draw[thick,dashed] (3,0) -- (3,0.5)node[left]{$2.A_4$};
            \draw[fill=white, draw=black, line width=0.025cm] (0,0) circle (0.02) node[below]{$2.A_4$};
            \draw[fill=white, draw=black, line width=0.025cm] (1,0) circle (0.02) node[below]{$C_{12}.2$};
            \draw[fill=white, draw=black, line width=0.025cm] (2,0) circle (0.02) node[below]{$C_{12}.2$};
            \draw[fill=white, draw=black, line width=0.025cm] (3,0) circle (0.02) node[below]{$2.A_4$};
        \end{tikzpicture}
    \end{center}
    
    Each $2.A_4$ endpoint has two remaining edges with stabilizer $2.A_4$. Since a $2.A_4$ fixed tree is finite, these cannot both join the two endpoints. Exactly one such edge would force the fixed subtree of $2.A_4$ to have length at least $3$, contradicting to $r_2=2$. Thus they lead to additional vertices with stabilizers of order at most $120$. Summing over vertices of the same type gives
    \[
    \mu(\Gamma\backslash G)\geq
    \frac1{24}+\frac1{24}+\frac1{120}+\frac1{120}
    =\frac1{10},
    \]
    equal to the stated upper bound. The only resulting graph of groups with such a covolume is illustrated in Example \ref{exm:y}. When lifting this graph of groups to the tree, the uniqueness of the $C_{12}$-subgroup in $C_{12}.2$ results in a bi-infinite path fixed by the same order-$12$ element, so Theorem \ref{thm:regular_width} implies $i\in K$. Since $q=4$ and $r_2=2$, this happens only when $K = \mathbb{Q}_2(\xi_3,\xi_4) = \mathbb{Q}_2(\xi_{12})$.

    \textbf{Case (5): $q = 2$.} In this case, vertex stabilizers are $Q_8$, $C_6.2$, $C_{12}$, $C_{12}.2$, and $2.A_4$. If there are no vertex subgroups isomorphic to $Q_8$, then in the absence of the exceptional star-orbit configuration for this $Q_8\cong C_{4(q+1)/3}.2$ or $C_{4(q-1)}.2$, the same argument as in Case (1) still holds, and the only additional case is Example \ref{exm:5}. Here, the elements of order $4$ in $Q_8$ act on the star for $C_{12}.2$ with order $2$, so $T_{Q_8}$ has diameter $r_2 = 2$. This only holds when $K = \mathbb{Q}_2(\sqrt{3})$.

    Suppose there is a $Q_8$ vertex. Two edges between two $Q_8$-vertices are not possible: two $Q_8$ edges will cause an infinite $Q_8$-fixed tree, while a $C_4$ edge and a $Q_8$ edge contradict themselves on whether the two $Q_8$ vertex groups are distinct. 

    Suppose first that every vertex has distance $\leq 1$ from a certain $Q_8$ vertex $v$. Without the prescribed double edges, the edge groups in the star of the $Q_8$ vertex, either $(Q_8,C_4)$ or $(Q_8,Q_8,Q_8)$, have vertex groups $C_6.2,C_{12}>C_4$, or $C_{12}.2,2.A_4>Q_8$ at the other endpoints. If the edge stabilizers are $(Q_8,Q_8,Q_8)$, at least two neighboring vertex groups are $C_{12}.2$. By Propositions \ref{prop:supergroup} and \ref{prop:dicyclic_fixed}, their distance equals the diameter of $T_{Q_8}$, hence $r_2=2$. Lemma \ref{lem:exceptional_tree_structure} then excludes $2.A_4$ at the third neighbor. If the edge stabilizers are $(Q_8,C_4)$, the $Q_8$ vertex is a leaf of $T_{Q_8}$, while its neighbor along the $Q_8$-edge is either a leaf or the center, by Proposition \ref{prop:supergroup} and Lemma \ref{lem:exceptional_tree_structure}. Thus $r_2\leq 2$.
    
    Now since $C_{12}$ and $C_{12}.2$ are exclusive to the case $K = \mathbb{Q}_2(\sqrt{3})$ with a smaller covolume, Example \ref{exm:6} is the only remaining possibility, where $r_2 = 2$. In this case, as the $C_4$ edge group is not normal in either of the vertex groups $2.A_4$ and $C_6.2$, the fixed subtree of the $C_4$ edge group turns out to be the $1$-neighborhood of the $Q_8$ edge. This requires that $K(i)/K$ is ramified, so $K\neq \mathbb{Q}_2(\sqrt{3})$ or $\mathbb{Q}_2(\sqrt{-1})$.

    If there is a vertex $w$ of distance $2$ from the $Q_8$-vertex $v$, the graph of groups will give the covolume bound
    \[
    \mu(\Gamma\backslash G) \geq \frac{1}{8} + \frac{1}{24} = \frac{1}{6},
    \]
    equal to the upper bound. This upper bound is attained only if the only vertex $w$ with distance $2$ from $v$ is either $2.A_4$ or $C_{12}.2$. In either case, the only outgoing edge from $w$ has a $Q_8$ group, and the vertex $u$ between $v$ and $w$ also has a $Q_8$ group. In addition, the uniqueness of $w$ with $d(w,v)\ge 2$ implies there are only two edges incident with $u$, so $\Gamma_{[v,u]}\cong C_4$. The other outgoing edge from $v$ then has stabilizer $Q_8$, and the covolume $\frac{1}{6}$ then implies the stabilizer at the other endpoint is either $2.A_4$ or $C_{12}.2$, as well. The $C_{12}.2$ case is impossible: otherwise the $Q_8$-fixed tree is an edge, so $r_2 = 1$ and $K = \mathbb{Q}_2$, contradicting to that $z_{12}\in K$. For the $2.A_4$ case described in Example \ref{exm:x}, the fixed tree of $Q_8$ implies $r_2 = 2$. Lifting the graph of groups to the covering tree, the normality of $C_4<Q_8$ gives a bi-infinite path fixed by the same order-$4$ element. Similarly to the explanation for Example \ref{exm:y}, we have $\xi_4\in K$, so $K = \mathbb{Q}_2(\xi_4)$.
\end{proof}
\begin{lem}[cf. {\cite[Lemma 7.1 and 7.2]{lubotzky1999lattices}}]\label{lem:order_8}
    Let $q=2,4,8$ or $16$, and let $\Gamma$ be a lattice of minimal covolume containing an element $g$ of order $8$. Then $\Gamma$ satisfies one of the following:
    \begin{itemize}
        \item $d_3\leq 1$, $r_3\leq 3$, and $K(i)/K$ is ramified; or
        \item $q\leq 4$ and there is a vertex $x\in\mathcal{V}(T)$ with $\Gamma_x\cong Q_{16}$; $r_3\leq 5$ if $q=4$, and $r_3\leq 7$ if $q=2$.
    \end{itemize}
\end{lem}
\begin{proof}
    Suppose first that $d_3\leq 1$. If $\xi_{2^3}\in K$ or $K(\xi_{2^3})/K$ is unramified, Theorem \ref{thm:regular_width} implies that $r_3 = 1$. However, in this case we have
    \[
    v_K(\xi_{2^3}-1) = \frac{e(K/\mathbb{Q}_2)}{2^{3-1}} = \frac{1}{2},
    \]
    contradicting the assumption. So $K(\xi_{2^3})/K$ is ramified, and Theorem \ref{thm:width_formula} implies that $r_3\leq 3$.

    Now suppose $d_3\geq 2$, so there exists a vertex $x\in \mathcal{V}(T_g) - (\partial T_g\cup \mathcal{V}(\mathfrak{C}_g))$. Take a vertex $x$ such that $d(x,\mathfrak{C}_g) = d_3 - 1$. We claim that $\Gamma_x$ does not contain an element of order $2^5$. Suppose the contrary. Proposition \ref{prop:not_boundary} (b) implies that $\Gamma_x$ is cyclic or dicyclic, so there is $h\in\Gamma_x$ with $h^4 = g$. By Theorem \ref{thm:regular_width} or \ref{thm:width_formula},
    \[
    (d_3-1) - d_5 = (d_3-d_5)-1 = (r_3-r_5)-1 = 3r_5-1>0,
    \] 
    a contradiction. 
    
    Since $g\in \Gamma_x$, by Proposition \ref{prop:not_boundary} (b), $\Gamma_x$ is a subgroup of $Q_{32}$. We will show that $\Gamma_x$ does not contain an element of order $2^4$ either: otherwise, there is $h\in\Gamma_x$ with $h^2 = g$. In particular, 
    \[
    1 = d_3-d_4 = r_3 - r_4 = r_4,
    \]
    so $r_3 = 2$ and $r_4 = 1$. Therefore,
    \[
    v_K(\xi_{2^4}-1) = \frac{e(K/\mathbb{Q}_2)}{2^{4-1}} = \frac{1}{2},
    \]
    meaning that $K(\xi_{2^4})/K$ is ramified. Theorem \ref{thm:width_formula} then implies that $d_4 = 0$ and $d_3 = 1$, a contradiction.

    Next suppose that $\Gamma_x$ is isomorphic to $C_8$. Since $\Gamma$ has minimal covolume, we have $q=2$. For a vertex $y\in \partial T_g$ adjacent to $x$, let $e = [x,y]$. Then $\Gamma_e\cong C_8$. Also, $\Gamma_y\ge \Gamma_e$ and $\lvert \Gamma_y\rvert \le (q+1)\lvert \Gamma_e\rvert$, implying that $\Gamma_y$ is cyclic or dicyclic; thus $\langle g\rangle$ is normal in $\Gamma_y$. Since $e$ is the only edge in $\st(y)$ fixed by $g$, $\Gamma_y$ fixes $e$, so $\Gamma_y \cong C_8$ as well. 
    
    Now there are two such vertices $y$ and $y'$, and they are not in the same $\Gamma$-orbit: if $\gamma\in \Gamma$ is such that $y = \gamma. y'$, then $\gamma\langle g\rangle\gamma^{-1} = \gamma \Gamma_{y'}\gamma^{-1} = \Gamma_{y} = \langle g\rangle$, so $\gamma$ preserves $T_g$; therefore $\gamma$ fixes $x$, contradicting the fact that $\Gamma_x = \langle g\rangle$. Therefore, $\mu(\Gamma\backslash G)\geq \frac{1}{\lvert \Gamma_{y}\rvert}+ \frac{1}{\lvert \Gamma_{y'}\rvert}\geq \frac{1}{4}$, contradicting the minimal covolume assumption.

    Therefore, when $d_3\geq 2$, the only possibility is that $\Gamma_x\cong Q_{16}$, so the minimal covolume assumption requires $q=2$ or $q=4$. For $q=2$, if $d_3\leq 3$, then Theorem \ref{thm:width_formula} implies $r_3\leq 7$. If $d_3\geq 4$, taking a vertex $u\in \mathcal{V}(T_g)$, with $d(u,\mathfrak{C}_g) = d_3-3$, we claim that $\Gamma_u$ contains an element of order $16$: otherwise, Proposition \ref{prop:not_boundary}(b) implies that $\Gamma_u\le Q_{16}$. Consider the set of vertices
    \[
    \Sigma = \{v\in \mathcal{V}(T_g): [v,\mathfrak{C}_g]\ni u,\ d(u,v) = 2\},
    \]
    Then $\lvert \Sigma\rvert = 4$, and $\Sigma \subset \mathcal{V}(T_g) - \partial T_g$. For $v,v'\in\Sigma$, arguments similar to those above imply that $\langle g\rangle$ is normal in $\Gamma_v$, so $\Gamma_v$ fixes $u$ and $\Gamma_v\le \Gamma_u$; and each $\gamma\in\Gamma$ with $v = \gamma.v'$ satisfies $\gamma\in \Gamma_u$. Therefore, we have
    \[
    \mu(\Gamma\backslash G)\geq \sum_{v\in \Gamma_u\backslash \Sigma}\frac{1}{\lvert \Gamma_v\rvert} = \frac{\lvert\Sigma\rvert}{\lvert\Gamma_u\rvert}\geq \frac{1}{4},
    \]
    contradicting the minimal covolume assumption. The existence of an element of order $16$ in $\Gamma_u$ then implies $d_4\geq d_3-3$, so $r_4 = d_3-d_4\leq 3$ and $r_3\leq 6<7$. For $q=4$, a similar argument proves that $r_3\leq 5$.
\end{proof}
\begin{thm}[cf. {\cite[Theorem 7.3, 7.4, and 7.5]{lubotzky1999lattices}}]\label{thm:order_8}
    Let $p = 2$ and let $\Gamma$ be a lattice of minimal covolume containing an element $g$ of order $8$. Then $q=2$ or $4$, and one of the following holds:

    \begin{enumerate}
        \item $K = \mathbb{Q}_2(z_8)$, $T//\Gamma$ is isomorphic to
        \begin{center}
            \begin{tikzpicture}[scale=2.5,every node/.style={font=\footnotesize}]
                \draw[thick] (0,0) --node[above]{$Q_{16}$} (1,0);
                \fill[] (0,0) circle(0.025) node[left]{$2.S_4$};
                \fill[] (1,0) circle(0.025) node[right]{$2.S_4$};
            \end{tikzpicture}
        \end{center}
        and $\mu(G/\Gamma) = \frac{1}{48}$.
        \item $K = \mathbb{Q}_2(z_8,z_{12})$, $T//\Gamma$ is isomorphic to
        \begin{center}
            \begin{tikzpicture}[scale=2.5,every node/.style={font=\footnotesize}]
                \draw[thick] (0,0) --node[above]{$Q_{16}$} (1,0) --node[above]{$Q_{16}$} (2,0);
                \draw[thick] (1,0) --node[left]{$Q_{16}$} (1,1);
                \fill[] (0,0) circle(0.025) node[below]{$2.S_4$};
                \fill[] (1,0) circle(0.025) node[below]{$Q_{16}$};
                \fill[] (2,0) circle(0.025) node[below]{$2.S_4$};
                \fill[] (1,1) circle(0.025) node[above]{$C_{24}.2$};
            \end{tikzpicture}
        \end{center}
        and $\mu(G/\Gamma) = \frac{1}{16}$.
        \item $K = \mathbb{Q}_2(z_{16})$, $T//\Gamma$ is isomorphic to
        \begin{center}
            \begin{tikzpicture}[scale=2.5,every node/.style={font=\footnotesize}]
                \draw[thick] (0,0) --node[above]{$Q_{16}$} (1,0) --node[above]{$Q_{32}$} (2,0) --node[above]{$Q_{16}$} (3,0);
                \fill[] (0,0) circle(0.025) node[below]{$2.S_4$};
                \fill[] (1,0) circle(0.025) node[below]{$Q_{32}$};
                \fill[] (2,0) circle(0.025) node[below]{$Q_{32}$};
                \fill[] (3,0) circle(0.025) node[below]{$2.S_4$};
            \end{tikzpicture}
        \end{center}
        and $\mu(G/\Gamma) = \frac{5}{96}$.
        \item $K = \mathbb{Q}_2(\xi_8)$, $T//\Gamma$ is isomorphic to
        \begin{center}
            \begin{tikzpicture}[scale=2.5,every node/.style={font=\footnotesize}]
                \draw[thick] (0,0) --node[above]{$Q_{16}$} (1,0) --node[above]{$Q_{16}$} (2,0) --node[above]{$C_8$} (3,0)--node[above]{$Q_{16}$} (4,0)--node[above]{$Q_{16}$} (5,0);
                \draw[thick] (1,0) --node[left]{$Q_{16}$}(1,1);
                \draw[thick] (4,0) --node[left]{$Q_{16}$}(4,1);
                \fill[] (0,0) circle(0.025) node[below]{$2.S_4$};
                \fill[] (1,0) circle(0.025) node[below]{$Q_{16}$};
                \fill[] (2,0) circle(0.025) node[below]{$Q_{16}$};
                \fill[] (3,0) circle(0.025) node[below]{$Q_{16}$};
                \fill[] (4,0) circle(0.025) node[below]{$Q_{16}$};
                \fill[] (5,0) circle(0.025) node[below]{$2.S_4$};
                \fill[] (1,1) circle(0.025) node[above]{$2.S_4$};
                \fill[] (4,1) circle(0.025) node[above]{$2.S_4$};
            \end{tikzpicture}
        \end{center}
        and $\mu(G/\Gamma) = \frac{1}{6}$.
        \item $K = \mathbb{Q}_2(z_5,z_8)$, $T//\Gamma$ is isomorphic to
        \begin{center}
            \begin{tikzpicture}[scale=2.5,every node/.style={font=\footnotesize}]
                \draw[thick] (0,0) --node[above]{$2.A_4$} (1,0) --node[above]{$Q_{16}$} (2,0) --node[above]{$2.A_4$} (3,0);
                \fill[] (0,0) circle(0.025) node[below]{$2.A_5$};
                \fill[] (1,0) circle(0.025) node[below]{$2.S_4$};
                \fill[] (2,0) circle(0.025) node[below]{$2.S_4$};
                \fill[] (3,0) circle(0.025) node[below]{$2.A_5$};
            \end{tikzpicture}
        \end{center}
        and $\mu(G/\Gamma) = \frac{7}{240}$.
    \end{enumerate}
\end{thm}
\begin{proof}
    By Example \ref{exm:reg}, there exists a lattice with covolume $\frac{1}{2(q+1)}$. If the lattice $\Gamma$ contains an element $g$ of order $8$, then for any $x\in\partial T_g$, Proposition \ref{prop:120_bound} implies that $\lvert \Gamma_x\rvert\leq 48$. This shows that $\Gamma$ is a lattice of minimal covolume only if $q\in \{2,2^2, 2^3, 2^4\}$.

    If $q = 2^3$ or $2^4$, then for $x\in\partial T_g$ described above, $\Gamma_x\cong C_8$, $Q_{16}$ or $2.S_4$, and only $2.S_4$ satisfies the minimal covolume requirement. By Proposition \ref{prop:edge_stab}, there is an edge $e = [x,y]$ incident with $x$, such that $\Gamma_e\cong Q_8$. Since $2.S_4$ acts on $\st(x)^\circ$ with a kernel isomorphic to $Q_8$, this kernel equals $\Gamma_e$ and is normal in $\Gamma_x$. By Lemma \ref{lem:order_8}, $e(K/\mathbb{Q}_2) = 2r_3\le 6$; since $v_2(z_{2^5}) = \frac{1}{8}<\frac{1}{6}$, $K$ does not contain elements of order $2^5$. The lemma also requires that $K(i)/K$ is ramified, so by Corollary \ref{cor:ramify_no_odd}, the vertex group $\Gamma_y$ can only be either exceptional, $Q_{16}$ or $Q_{32}$. 
    
    Here, if $q = 2^4$, minimal covolume requires $\lvert \Gamma_y\rvert \geq 34$, while Proposition \ref{prop:edge_stab} implies $\Gamma_y\not\cong 2.A_5$. Thus we have $\Gamma_y\cong 2.S_4$. However, $\Gamma_e$ is then also normal in $\Gamma_y$, while Proposition \ref{prop:supergroup_exceptional} only allows one vertex group isomorphic to $2.S_4$ in $T_{\Gamma_e}$ that contains $\Gamma_e$ as a normal subgroup, a contradiction. If $q = 2^3$, then similarly $\Gamma_y\not\cong 2.S_4$, while $Q_{32}$ does not allow an edge stabilizer isomorphic to $Q_8$. Since $2\nmid f$, Proposition \ref{prop:supergroup_exceptional} implies that only the center of $T_{\Gamma_e}$ could have vertex group $2.A_4$ or $2.S_4$ in which $\Gamma_e$ is normal, so $\Gamma_y\cong 2.A_4$ is also impossible. Therefore we must have $q=2$ or $q = 2^2$, which are treated in the propositions below.
\end{proof}
\begin{prop}\label{prop:q_is_4}
    Under the hypotheses of Theorem \ref{thm:order_8}, if $q=4$, then statement (5) holds.
\end{prop}
\begin{proof}
    By Lemma \ref{lem:order_8}, $r_3\leq 5$, so we distinguish cases according to $r_3$. In particular, if $r_3>3$, a vertex group $Q_{16}$ is required; in other cases, Proposition \ref{prop:120_bound}, the order-$8$ restriction and the covolume restriction also require a vertex group of either $Q_{16}$ or $2.S_4$.

    \textbf{Case (i): $r_3 = 1$ ($r_2 = 2$).} This case occurs only when $K = \mathbb{Q}_2(z_8,z_5)$. In this case, $K(i)/K$ is quadratic and ramified, so $d_3 = 0$ and $d_2 = 1$; since there are no elements of order $16$, $12$ or $20$, only $Q_{16}$, $C_6.2$, $C_{10}.2$, $2.A_4$, $2.S_4$ and $2.A_5$ are allowed as vertex groups. Since a vertex group of order at most $\lvert 2.S_4\rvert = 48$ exists and Proposition \ref{prop:no_8} implies $\lambda'(K) = \frac{1}{10}$, we have $\lambda(K)\ge \frac{1}{48}$. This lower bound is not attained: a $2.S_4$ vertex group has an edge group of $Q_{16}$, which has to connect to another $2.S_4$ vertex group by the list of possible vertex groups and the covolume bound. Then the other $2.A_4$ edge must also connect to this $2.S_4$ vertex. However, the fixed subtree of $2.A_4$ cannot contain two vertices with stabilizer $2.S_4$, a contradiction. Therefore, the next possible covolume is $\frac{1}{48}+\frac{1}{120} = \frac{7}{240}$, which is attained in Example \ref{exm:b}.

    \textbf{Case (ii): $r_3 \in \{3,5\}$.} Since $r_3$ is odd, $\slg(2,K)$ does not have elements of order $2^4 = 16$; since $K(i)\supset \mathbb{Q}_2(z_8,\xi_4)$ and $e(\mathbb{Q}_2(z_8,\xi_4)/\mathbb{Q}_2) = 4$, we see that $K(i)/K$ is ramified. By Corollary \ref{cor:ramify_no_odd}, $\slg(2,K)$ does not contain $C_{12}$ or $C_{20}$ either. Therefore, the vertex group candidates are the same as in Case (i).

    Suppose there is a vertex group isomorphic to $2.S_4$, so the edge groups are $Q_{16}$ and $2.A_4$. On the $Q_{16}$ side, since $r_3>1$, Proposition \ref{prop:supergroup_exceptional} implies that the other endpoint cannot have stabilizer $2.S_4$. Therefore, this vertex group has to be $Q_{16}$. On the $2.A_4$ side, a $2.A_4$ vertex group at the other endpoint would give a covolume of at least $\frac{1}{16}+\frac{1}{24} = \frac{5}{48}>\frac{1}{10}$; a $2.S_4$ vertex group is excluded as in Case (i). Since $r_3>1$, Corollary \ref{cor:A4_attach} implies $2.A_5$ is also impossible.

    Therefore one has to start with a vertex group of $Q_{16}$, without using $2.S_4$ as a vertex group. The $Q_{16}$ vertex group has a $Q_{16}$ edge, which has to connect to another $Q_{16}$ vertex group. The fact that $r_3>1$ again forces one of the two vertices to be an interior vertex of the fixed subtree of this $Q_{16}$ group. This vertex is incident with five $Q_{16}$ edges in the quotient graph of groups. Without using $2.S_4$, all of these edges should connect to $Q_{16}$ vertices. The finiteness of $T_{Q_{16}}$ prevents any double $Q_{16}$ edges between $Q_{16}$ vertices, so the covolume is at least $\frac{5}{16}>\frac{1}{10}$.

    \textbf{Case (iii): $r_3 = 2$.} For the cases where $K\neq \mathbb{Q}_2(z_{16},z_5)$, $\mathbb{Q}_2(z_8,z_{20})$ or $\mathbb{Q}_2(z_5,z_8,z_{12})$, $\slg(2,K)$ does not contain elements of orders $16$, $12$ or $20$, and the vertex group candidates are still the same as in Case (i); this is excluded similarly to Case (ii). We then focus on the three cases above.

    \textit{Case (iii.a): $K = \mathbb{Q}_2(z_{16},z_5)$.} Here $K(i)/K$ is ramified, so $d_3 = 1$ and $d_2 = 3$; $C_{16}$ and $Q_{32}$ are also candidate vertex groups. If there is a $2.S_4$ vertex $v_1$, we have seen that the $2.A_4$ edge can only connect to a $2.A_4$ vertex $v_0$. For the $Q_{16}$ edge, we have seen that $Q_{16}$ or $2.S_4$ vertices are not possible, so it connects to a $Q_{32}$ vertex $v_2$. The $Q_{32}$ edge incident with this vertex connects to another $Q_{32}$ vertex $v_3$. If a $C_{16}$ or $C_{8}$ edge is incident with $v_3$, then its other endpoint is a new vertex with cyclic or dicyclic stabilizer of order at most $32$; the covolume is at least $\frac{1}{24} + \frac{1}{32} + \frac{1}{32} > \frac{1}{10}$. If two more $Q_{16}$ edges are incident with $v_3$, at most one of them can connect to $v_2$, and the other has to connect to a $2.S_4$ vertex $v_4$ by the covolume bound. Now the $2.A_4$ edge of $v_4$ similarly connects to a $2.A_4$ vertex $v_5$; since $r_3=2$ and this $2.A_4$ vertex is adjacent to a $2.S_4$ vertex, it is in the interior of the fixed subtree of its stabilizer; so there are $2.A_4$ and $Q_8$ edges incident with $v_5$; in particular, the $Q_8$ edge cannot be incident with $v_0$. We have the following configuration:
    \begin{center}
        \begin{tikzpicture}[scale=2.5,every node/.style={font=\footnotesize}]
            \draw[thick,dashed] (-0.5,0) -- (0,0) --node[above]{$2.A_4$} (1,0) --node[above]{$Q_{16}$} (2,0) --node[above]{$Q_{32}$} (3,0) --node[above]{$Q_{16}$} (4,0)--node[above]{$2.A_4$} (5,0) -- (5.5,0);
            \draw[thick,dashed] (2,0) -- (2,0.5);
            \draw[thick,dashed] (3,0) -- (3,0.5);
            \draw[fill=white, draw=black, line width=0.025cm] (0,0) circle(0.025) node[below]{$2.A_4$};
            \draw[fill=white, draw=black, line width=0.025cm] (1,0) circle(0.025) node[below]{$2.S_4$};
            \draw[fill=white, draw=black, line width=0.025cm] (2,0) circle(0.025) node[below]{$Q_{32}$};
            \draw[fill=white, draw=black, line width=0.025cm] (3,0) circle(0.025) node[below]{$Q_{32}$};
            \draw[fill=white, draw=black, line width=0.025cm] (4,0) circle(0.025) node[below]{$2.S_4$};
            \draw[fill=white, draw=black, line width=0.025cm] (5,0) circle(0.025) node[below]{$2.A_4$};
        \end{tikzpicture}
    \end{center}
    
    The remaining $2.A_4$ edge at $v_5$ leads to a new vertex with stabilizer of order at most $120$, hence $\mu(\Gamma\backslash G) \geq \frac1{24}+\frac1{32}+\frac1{48}+\frac1{120} =\frac{49}{480}>\frac1{10}$.

    If there is a $Q_{16}$ vertex $v_1$ and no $2.S_4$ vertices, then the $Q_{16}$ edge connects to either a $Q_{32}$ vertex or a $Q_{16}$ vertex $v_2$. If $v_2$ has stabilizer $Q_{16}$, then either $v_1$ or $v_2$ is an interior vertex of the fixed subtree of this $Q_{16}$ group. The vertex has five $Q_{16}$ edges incident with it; the other four $Q_{16}$ edges connect to either $Q_{16}$ or $Q_{32}$ groups, and each $Q_{32}$ vertex is incident with at most two $Q_{16}$ edges. Therefore, the covolume is at least $\frac{1}{16}+ \frac{2}{32}>\frac{1}{10}$. If $v_2$ has vertex group $Q_{32}$, the other two edge groups at $v_2$ are $Q_{16}$ and $Q_{32}$. The $Q_{32}$ edge has to connect to another $Q_{32}$ vertex $v_3$. The $Q_{16}$ edge cannot connect to $v_1$, otherwise the uniqueness of $Q_{32}$ in $\slg(2,K)$ containing $Q_{16}$ contradicts the finiteness of the $Q_{16}$ fixed subtree. Nor can it connect to $v_3$, otherwise the uniqueness of $Q_{32}$ containing $Q_{16}$ forces the vertex groups at $v_2$ and $v_3$ to coincide, and a $Q_{16}$ edge in between is not possible. Therefore this $Q_{16}$ edge must connect to a new vertex, giving a covolume greater than $\frac{1}{10}$ again.

    \textit{Case (iii.b): $K = \mathbb{Q}_2(z_{8},z_{20})$.} Here $K(i)/K$ is unramified, so $d_3 = 2$ and $d_2 = 4$, and a vertex stabilizer isomorphic to $Q_{16}$ exists. Additional candidate vertex groups are $C_{20}$, $C_{20}.2$, $C_{40}$, and $C_{40}.2$. In this case, the center of a $Q_{16}$ fixed subtree cannot admit a larger vertex group such as $2.S_4$, so there must be a $Q_{16}$ vertex with five $Q_{16}$ edges. At most two edges are incident with $2.S_4$ vertices, while at most one is incident with a $C_{40}.2$ vertex, so at least two edges are incident with $Q_{16}$ vertices, and these vertices must be different. The covolume is now greater than $\frac{1}{10}$.
    
    \textit{Case (iii.c): $K = \mathbb{Q}_2(z_{8},z_{5},z_{12})$.} Here $i\in K$, so $d_3 = 2$ and $d_2 = 4$, and a vertex stabilizer isomorphic to $Q_{16}$ exists. Similarly to Case (iii.b), this is not possible either.
    
    \textbf{Case (iv): $r_3 = 4$.} In this case, the required vertex group isomorphic to $Q_{16}$ has a fixed subtree of diameter $4$. Since the $Q_{64}$ supergroup, if it exists, fixes only an edge of $T_{Q_{16}}$ that meets $\partial T_{Q_{16}}$, the center $v$ of $T_{Q_{16}}$ has vertex stabilizer $Q_{16}$ or $Q_{32}$. If $\Gamma_v\cong Q_{16}$, then there are five $Q_{16}$ edges in the quotient graph of groups, with $Q_{16}$ groups at the other endpoints, and these endpoints are distinct. If $\Gamma_v\cong Q_{32}$, since $v\in\partial T_{Q_{32}}$, there are still two $Q_{16}$ edges incident with it. These edges cannot connect to the same $Q_{16}$ vertex either: if the tree vertices $x,y$ have the same image in the quotient graph of groups, say $g.x = y$ for some $g\in\Gamma$, then $\Gamma_x=\Gamma_y\cong Q_{16}$ implies that $g$ normalizes this $Q_{16}$ group, so $g$ fixes the center $v$, contradicting the fact that $[v,x]$, $[v,y]$ are in different $\Gamma_v$-orbits. This gives the covolume bound $\frac{2}{16}>\frac{1}{10}$.
\end{proof}
\begin{prop}\label{prop:q_is_2}
    Under the hypotheses of Theorem \ref{thm:order_8}, if $q=2$, then one of statements (1), (2), (3) and (4) holds.
\end{prop}
\begin{proof}
    Similarly, $r_3\le 7$, either a $Q_{16}$ or a $2.S_4$ vertex group is required, and if $r_3>3$, a $Q_{16}$ vertex group is required.

    \textbf{Case (i): $r_3=1$.} The only case is $K = \mathbb{Q}_2(z_8)$; since a group of order at most $\lvert 2.S_4\rvert=48$ is required, the construction in Example \ref{exm:a} gives the unique graph of groups for lattices of minimal covolume in $\slg(2,K)$; the covolume is $\frac{1}{48}$.

    \textbf{Case (ii): $r_3\in\{3,5,7\}$.} If $r_3=3$ and there is a $2.S_4$ vertex $v_1$, then the $Q_{16}$ edge connects to a $Q_{16}$ vertex $v_2$; in particular, this $v_2$ is an interior vertex of the fixed subtree of this $Q_{16}$ group and is incident with two additional $Q_{16}$ edges. One of them connects to another interior vertex $v_3$, i.e., $[v_2v_3]$ is the central edge of $T_{Q_{16}}$. If all ($Q_{16}$-) edges incident with $v_3$ connect to $Q_{16}$ vertices, they are distinct vertices, and the covolume is at least $\frac{3}{16}>\frac{1}{6}$. Therefore, one of the vertices adjacent to $v_3$ has group $2.S_4$. Thus the graph of groups contains the following configuration:
    \begin{center}
        \begin{tikzpicture}[scale=2.5,every node/.style={font=\footnotesize}]
            \draw[thick,dashed] (-0.5,0) -- (0,0) --node[above]{$Q_{16}$} (1,0) --node[above]{$Q_{16}$} (2,0)--node[above]{$Q_{16}$} (3,0) -- (3.5,0);
            \draw[thick, dashed] (1,0) --node[left]{$Q_{16}$} (1,1);
            \draw[thick, dashed] (2,0) --node[left]{$Q_{16}$} (2,1);
            \draw[fill=white, draw=black, line width=0.025cm] (0,0) circle(0.025) node[below]{$Q_{16}$};
            \draw[fill=white, draw=black, line width=0.025cm] (1,0) circle(0.025) node[below]{$Q_{16}$};
            \draw[fill=white, draw=black, line width=0.025cm] (2,0) circle(0.025) node[below]{$Q_{16}$};
            \draw[fill=white, draw=black, line width=0.025cm] (3,0) circle(0.025) node[below]{$Q_{16}$};
            \draw[fill=white, draw=black, line width=0.025cm] (1,1) circle(0.025) node[above]{$2.S_4$};
            \draw[fill=white, draw=black, line width=0.025cm] (2,1) circle(0.025) node[above]{$2.S_4$};
        \end{tikzpicture}
    \end{center}
    The $Q_{16}$ vertices on the two sides are boundary vertices and are contained in additional $C_8$ or $Q_8$ edges. If they are incident to additional vertices, such a group has order at most $3\cdot 8 = 24$, making the covolume be at least $\frac{1}{24}+\frac{2}{16}+\frac{1}{48} = \frac{3}{16}>\frac{1}{6}$. Now suppose these $Q_{16}$ vertices are connected together by such an edge with group $C_8$ or $Q_8$. In this setting, $Q_{16}$ contains a unique $C_8$ subgroup, and Theorem \ref{thm:width_formula} and Proposition \ref{prop:dicyclic_fixed} imply the groups $C_8$ and $Q_{16}$ have the same fixed subtree: both are the $1$-neighborhood of an edge. This implies that an additional $C_8$ edge beyond $\partial T_{Q_{16}}$ is not possible. Now suppose they are connected by a $Q_8$ edge. Lifting this $Q_8$ edge to the tree, we note that the $Q_{16}$ supergroups fixing its two endpoints are different. Therefore, the graph of groups implies that these two $T_{Q_{16}}$ are at distance $1$ from each other, contradicting Proposition \ref{prop:dicyclic_fixed}, which requires them to meet at the center of $T_{Q_{8}}$.

    The remaining cases are either $r_3=3$ without $2.S_4$ vertices, or $r_3\geq 5$ and with a $Q_{16}$ vertex. In both cases, if we begin with the central edge of this $T_{Q_{16}}$, both endpoints $u,v$ are interior vertices. If $r_3 = 3$, we have forced these vertex groups to be $Q_{16}$; if $r_3 = 5$ or $7$, these vertices are still interior vertices, so the vertex groups should also be $Q_{16}$. In both cases, the four $Q_{16}$ edges incident with either of them do not share vertices other than $u,v$, which gives a covolume of at least $\frac{3}{16}>\frac{1}{6}$.
    
    \textbf{Case (iii): $r_3 = 2$.} If $K \neq \mathbb{Q}_2(z_8,z_{12})$ or $\mathbb{Q}_2(z_{16})$, then $\slg(2,K)$ has the same vertex group candidates as above. Whether starting with a $Q_{16}$ or a $2.S_4$ vertex, the existence of a $Q_{16}$ edge will force a $Q_{16}$ vertex $v$ at the fixed subtree center. Among the three edges incident with $v$, at most two of them could connect to $2.S_4$ vertices, so at least one of them connects to a vertex, denoted by $w$, with group $Q_{16}$. The other edge $[w,u]$ incident with $w$ has a $Q_8$ or $C_8$ group.
    
    If the vertex group of $u$ is $2.A_4$, then $u$ is the center of the fixed subtree of $\Gamma_{[w,u]}\cong Q_8$; hence $u$ is a boundary vertex of $T_{\Gamma_w}$. Since $w$ is also a boundary vertex of this ball of radius $1$, we have $d(w,u)\in\{0,2\}$, contradicting $d(w,u)=1$. The covolume bound then forces $\Gamma_u\cong Q_{16}$, and the other edge incident with $u$ connects to the center of $T_{\Gamma_u}$, corresponding to another $Q_{16}$ group. These $Q_{16}$ vertices belong to distinct orbits, and the two remaining $Q_{16}$ edges at the center of $T_{\Gamma_u}$ lead to two new vertices with stabilizers of order at most $48$. Together with $v$ and $u$, these give a covolume of at least $\frac{2}{48}+ \frac{2}{16} = \frac{1}{6}$, equal to the upper bound $\frac{1}{2(q+1)}$. For a similar reason as in Case (ii), the edge stabilizer between the two $Q_{16}$-fixed trees in this configuration cannot be $Q_8$. Therefore $\Gamma_{[w,u]}\cong C_8$, as Example \ref{exm:z} described. When lifting the graph of groups into a tree, the uniqueness of $C_8$-subgroup in $Q_{16}$ results in a bi-infinite path fixed by the same order-$8$ element, so Theorem \ref{thm:regular_width} requires $\xi_8\in K$, or $K = \mathbb{Q}_2(\xi_8)$.

    \textit{Case (iii.a): $K = \mathbb{Q}_2(z_{16})$.} A lattice with covolume $\frac{5}{96}$ is constructed in Example \ref{exm:d}, and we claim this is the minimal covolume: since $\frac{5}{96}<\frac{1}{16}$, we must start with a $2.S_4$ vertex $v$ and exclude $Q_{16}$ vertices. As $r_3=2$, the $Q_{16}$ edge incident with $v$ cannot connect to another $2.S_4$ vertex, so $Q_{32}$ is the only option for the stabilizer of the other endpoint $w$, apart from $Q_{16}$. The $Q_{32}$ edge incident with $w$ connects to another $Q_{32}$ group, giving $\lambda(K)\geq \frac{1}{48}+\frac{1}{32} = \frac{5}{96}$. It is clear that Example \ref{exm:d} gives the unique graph of groups for lattices of this minimal covolume in $\slg(2,K)$.  
    
    \textit{Case (iii.b): $K = \mathbb{Q}_2(z_8,z_{12})$.} In this case, $K(i)/K$ is unramified, and $d_3 = 2$. Lemma \ref{lem:order_8} forces a vertex group of $Q_{16}$, so the covolume is at least $\frac{1}{16}$. Propositions \ref{prop:supergroup} and \ref{prop:supergroup_exceptional} imply that $\partial T_{Q_{16}}$ allows at most two vertex groups isomorphic to $2.S_4$, and at most one dicyclic vertex group containing order-$3$ elements. Therefore, as the unique candidate with covolume $\frac{1}{16}$, the lattice with minimal covolume is constructed in Example \ref{exm:c}.
    
    \textbf{Case (iv): $r_3 = 6$.} Since $r_3>3$, there exists a vertex with group $Q_{16}$. $G$ contains no elements of order $32$, so the stabilizer of the center of this $T_{Q_{16}}$ is either $Q_{16}$ or $Q_{32}$. If $\Gamma_v \cong Q_{16}$, then $\Gamma$ does not contain the $Q_{32}$-supergroup of $\Gamma_v$, and all three edges incident with $v$ connect to $Q_{16}$ groups, making the covolume at least $\frac{3}{16}>\frac{1}{6}$. If $\Gamma_v\cong Q_{32}$, then there is an edge $[v,w]$, such that the stabilizers of both $[v,w]$ and the vertex $w$ coincide with the given $Q_{16}$ subgroup. Since $r_3=6$, $\st(w)$ is away from the boundary of this $Q_{16}$ tree, so the other two vertices in $T//\Gamma$ adjacent to $w$ also have $Q_{16}$ groups.

    On the other side, the $Q_{32}$-edge incident with $v$ connects to a $Q_{32}$ vertex $u$, which is an interior vertex of its fixed subtree. Therefore the other two edges incident with $u$ have stabilizer $Q_{32}$; one of them is an inner edge of $T_{Q_{32}}$, forcing another $Q_{32}$ vertex. The covolume is then at least $\frac{2}{16}+\frac{2}{32}>\frac{1}{6}$.
    
    \textbf{Case (v): $r_3 = 4$.} Similar to Case (iv), there is a vertex with group $Q_{16}$; similar to Case (iv) of Proposition \ref{prop:q_is_4}, the potentially existing $Q_{64}$ supergroup does not fix any inner edge of $T_{Q_{16}}$. The stabilizer of the center $v$ of $T_{Q_{16}}$ is then either $Q_{16}$ or $Q_{32}$, while the $Q_{16}$ case is excluded as in Case (iv).

    Now suppose that $\Gamma_v\cong Q_{32}$. In this case, we still have an edge $[v,w]$ such that $\Gamma_{[v,w]} = \Gamma_w\cong Q_{16}$. If a boundary edge incident with $w$ has its other endpoint at a $C_{24}.2$ vertex, then after lifting, there will be two such vertex groups in $T_{Q_{16}}$, contrary to Proposition \ref{prop:supergroup}. Since $r_3=4$, at most one of them can connect to a $2.S_4$ group. 
    
    Therefore, a vertex $u$ adjacent to $w$ admits a stabilizer isomorphic to $Q_{16}$. Notice that the group $\Gamma_u$ contains a $Q_8$-subgroup whose fixed subtree has center at a neighbor of $w$ other than $v$: an element of $\Gamma_v\setminus\Gamma_w$ exchanges the two $Q_8$ subgroups of $\Gamma_w$, and hence their fixed-tree centers. These centers lie at distance $2$ from $v$ in distinct directions. Since the third direction at $v$ is fixed by $\Gamma_v$, the two centers lie in the two directions exchanged by $\Gamma_v$, one of which passes through $w$. This subgroup fixes the star of $u$; thus the outgoing edge from $u$ has group $Q_8$. We now have the following configuration:
    \begin{center}
        \begin{tikzpicture}[scale=2.5,every node/.style={font=\footnotesize}]
            \draw[thick,dashed] (-0.5,0) node[above]{$Q_8$} -- (0,0) --node[above]{$Q_{16}$} (1,0) --node[above]{$Q_{16}$} (2,0) -- (2.5,0)node[above]{$Q_{32}$};
            \draw[thick,dashed] (1,0) -- (1,0.5) node[left]{$Q_{16}$};
            \draw[fill=white, draw=black, line width=0.025cm] (0,0) circle(0.025) node[below]{$Q_{16}$};
            \draw[fill=white, draw=black, line width=0.025cm] (1,0) circle(0.025) node[below]{$Q_{16}$};
            \draw[fill=white, draw=black, line width=0.025cm] (2,0) circle(0.025) node[below]{$Q_{32}$};
        \end{tikzpicture}
    \end{center}
    
    The other endpoint of this outgoing edge is at odd distance from the center of $T_{Q_8}$, so its stabilizer cannot be the $2.A_4$ supergroup. The same argument also prevents the other endpoint from having stabilizer $C_{12}.2$. Nor can this edge connect to a $Q_{16}$ vertex: two $Q_{16}$-supergroups of $Q_8$ have fixed subtrees meeting at the center of $T_{Q_8}$, without having an edge fixed by neither of them in between. Therefore, this endpoint has a vertex group $Q_8$, giving a covolume of at least $\frac{1}{8}+\frac{1}{16}>\frac{1}{6}$.
\end{proof}
Propositions \ref{prop:q_is_4} and \ref{prop:q_is_2} together complete the proof of Theorem \ref{thm:order_8}; the latter together with Proposition \ref{prop:no_8} proves our main Theorem \ref{thm:main}. This main theorem, with an emphasize of Examples \ref{exm:x}, \ref{exm:y} and \ref{exm:z}, implies Corollary \ref{cor:maincor}.

\textbf{Artificial Intelligence Declaration.} The mathematical ideas and proofs of this paper are from the author. Artificial Intelligence was only used for typesetting and language clean-up. We take the responsibility for the content.
\appendix
\section{Realization of the Width \texorpdfstring{$d_k$}{dk} in the Ramified Case}\label{append:a}
\begin{prop}\label{prop:actual_width}
    Under the hypotheses of Theorem \ref{thm:width_formula}, when $r_k$ is even, each integer between $r_k/2$ and $r_k-1$ is attained as $d_k$ for some choice of $K$. Moreover, let $\mathbb{F}_q$ be the residue field of $K$, $K_0\coloneqq \mathbb{Q}_p(\xi_{q-1},z_{2^k})$, and $\pi$ be a uniformizer, so $K = K_0(\pi)$, and $e(K/K_0) = r_k$. Put $w = w_{2^k} \coloneqq 2-z_{2^k}$, so $w$ is a uniformizer of $K_0$ and $v(w) = r_k v(\pi)$, and put $r=r_k$. Suppose
    \[
    f(x) = x^r - w(a_{r-1}x^{r-1}+\dots + a_1x+a_0)
    \]
    is an Eisenstein polynomial of $\pi$, where $a_{r-1}$,\dots, $a_0\in \mathcal{O}_{K_0}$ and $a_0\not\equiv 0\bmod w$. By abuse of notation, we use the same symbols for their residue classes in $\mathcal{O}_{K_0}/(w.\mathcal{O}_{K_0}) = \mathbb{F}_q$. Notice that $\mathbb{F}_q$ has characteristic $2$ and every element has a square root. Then the equations
    \[
    a_j = \sum_{i=1}^{\lfloor \log_2(r/(r-j))\rfloor}(a_{r-2^i(r-j)})^{1/2^i},\quad j=1,3,\dots, r-1,\ 2\nmid j,
    \]
    cannot all be satisfied. Let $l$ be the smallest index $j$ such that the equation does not hold. Then the width is $d_k = \frac{r_k+l-1}{2}$.
\end{prop}
\begin{lem}\label{lem:gap_formula}
    Suppose $\operatorname{char}(\mathbb{F}_q) = 2$, and let
    \[
    g(x) = a_{r-1}x^{r-1}+\dots + a_1x+a_0\in \mathbb{F}_q[x],
    \]
    be a fixed polynomial of degree at most $r-1$ with $a_0\neq 0$ and $2\mid r$, and let
    \[
    C(x) = c_{r/2-1}x^{r/2-1} + \dots + c_1x + c_0\in \mathbb{F}_q[x]
    \]
    be any polynomial of degree at most $r/2-1$. Then the smallest odd number such that the equation
    \[
    a_j = \sum_{i=1}^{\lfloor \log_2(r/(r-j))\rfloor}(a_{r-2^i(r-j)})^{1/2^i}
    \]
    does not hold, if any, is the largest $j\in\mathbb{N}$ such that
    \[
    x^j\mid (1+x^{r/2}C(x) + C^2(x)g(x))
    \]
    for some $C(x)\in \mathbb{F}_q[x]$.
\end{lem}
\begin{proof}
    Denote the key expression in the lemma by
    \[
    E(x) = 1+x^{r/2}C(x) + C^2(x)g(x).
    \]
    To simplify the equation, let
    \[
    V(x) = x^{-r/2}g(x)C(x),
    \]
    Then
    \[
    V^2(x) + V(x) = x^{-r}g(x)(1+E(x)).
    \]
    Express $V(x)$ as a Laurent series:
    \[
    V(x) = \sum_{i\leq r/2} v_{-i}x^{-i},\ v_{-i}\in \mathbb{F}_q.
    \]
    Since the terms of degree less than $j$ in $E(x)$ vanish, the $x^{-2l}$-term of the equation above is
    \[
    v_{-l}^2 + v_{-2l} = a_{r-2l}\implies v_{-l} = v_{-2l}^{1/2} + a_{r-2l}^{1/2},
    \]
    for each $l>(r-j)/2$. Since $2l>l$, iteration gives
    \[
    v_{-l} = \sum_{i= 1}^{\lfloor\log_2(r/l) \rfloor} a_{r - 2^il}^{1/2^i}.
    \]
    Now, for each odd number $l\leq j$, whether the $x^l$-term vanishes in $E(x)$ is equivalent to whether the $x^{l-r}$-coefficient equation above holds with the $E(x)$ component omitted. Since the $x^{l-r}$-term does not appear in $V^2(x)$, this is equivalent to whether $v_{l-r} = a_l$. Substituting the expression for $v_{-l}$ in terms of the coefficients $a_i$ yields the desired conclusion.
\end{proof}
\begin{proof}[Proof of Proposition \ref{prop:actual_width}]
    To show that each integer in this range occurs, we consider the Eisenstein polynomial
    \[
    x^r + wx^{r/2}+wx^j+w = 0,\quad j=1,3,\dots, r-1,\ 2\nmid j,
    \]
    and let $\pi$ be a root of this polynomial. Set $K_0 = \mathbb{Q}_p(\xi_{q-1},z_{2^k})$, and $K = K_0(\pi)$. Let $\beta = \pi^{r/2}$. Then
    \[
    v(\beta^2 + w\beta + w) = v(\pi^r + w\pi^{r/2} + w) = v(w\pi^j) = (r_k+j)v(\pi),
    \]
    so $d_k\geq \frac{r_k+j-1}{2}$.

    To show that there is no other $\beta'\in\mathcal{O}_K$ with $v(\beta'^2 + w\beta' + w)> (r_k+j)v(\pi)$, suppose $\beta' = \pi^{r/2} + y$, and $y\in\mathcal{O}_K$. Then
    \[
    \beta'^2 + w\beta' + w = (\pi^{r/2} + y)^2 + w(\pi^{r/2} + y) + w = w\pi^j + 2\pi^{r/2}y + y^2 + wy. 
    \]
    It holds that $v(2\pi^{r/2}y)>v(2y)\geq v(wy)$, so the term with the smallest valuation is among $w\pi^j$, $y^2$, and $wy$. Since $j$ is odd, $v(y^2)\neq r_k+j$. If the valuation is greater than $v(w\pi^j) = (r_k+j)v(\pi)$, we must have
    \[
    v(y^2)> (r_k+j)v(\pi),\text{ and }v(w\pi^j) = v(wy).
    \]
    Therefore $v(y) = jv(\pi)$, so $2j>r_k+j$ or $j>r_k$, a contradiction. This finishes the construction of $K$ such that $d_k = \frac{r_k+j-1}{2}$, which ranges over all integers between $\frac{r_k}{2}$ and $r_k-1$.

    To complete the classification, suppose the uniformizer $\pi\in K$ is a root of the Eisenstein polynomial
    \[
    f(x) = x^r - wg(x) = x^r - w(a_{r-1}x^{r-1}+\dots + a_1x+a_0),\ a_i\in \mathcal{O}_{K_0},\ a_0\in \mathcal{O}_{K_0}^\times.
    \]
    Also, suppose $v(\beta^2 + w\beta + w)$ attains its maximum at some $\beta$; from the reasoning above, we must have $v(\beta) = (r_k/2)v(\pi)$, so we assume
    \[
    \beta \equiv \pi^{r/2}C(\pi) = \pi^{r/2}(c_0 + c_1\pi + \dots c_{r/2-1}\pi^{r/2-1})\bmod \pi^r,\ c_0,\dots, c_{r/2-1}\in \mathcal{O}_{K_0}.
    \]
    By Lemma \ref{lem:width_ineq}, $v(\beta^2 + w\beta + w)\leq 2r_kv(\pi) = 2v(w)$, thus we can consider $\beta^2 + w\beta + w\bmod w^2$:
    \[
    \beta^2 + w\beta + w \equiv \pi^r C^2(\pi) + w\pi^{r/2}C(\pi) + w\equiv w(C^2(\pi)g(\pi) + \pi^{r/2}C(\pi) + 1)\bmod w^2.
    \]
    The expression modulo $w^2$ remains unchanged when each coefficient $a_i$ or $c_i$ changes by an element of $w\mathcal{O}_{K_0}$, so we can assume these coefficients are taken in $\mathcal{O}_{K_0}/w\mathcal{O}_{K_0} = \mathbb{F}_q$, in particular $a_0\neq 0$. Among different polynomials $C(x)$, the largest possible $v(\beta^2 + w\beta + w)$ corresponds to the smallest nonzero $x^j$ term in $1 + x^{r/2}C(x) + C^2(x)g(x)$. If some of the identities for $a_j$ fail, this is determined by Lemma \ref{lem:gap_formula}; suppose this is the $x^l$ term. Then
    \[
    \max_{\beta\in\mathcal{O}_K} v_K(\beta^2 + w\beta + w) = v(w\pi^l) = r_k+l,
    \]
    and $d_k = \frac{r_k+l-1}{2}$.
    
    On the other hand, if all the identities for $a_j$ held, the recursion in the proof of Lemma \ref{lem:gap_formula} would yield $C(x)$ with $E(x)\equiv0\pmod{x^r}$. Taking $\beta=\pi^{r/2}C(\pi)$, the congruence above would imply $v_K(\beta^2+w\beta+w)\geq2r$. But Lemmas \ref{lem:max_valuation} and \ref{lem:width_ineq} give
    \[
    \max_{\beta\in\mathcal O_K}v_K(\beta^2+w\beta+w)=2d_k+1\leq2r,
    \]
    which is impossible since $2d_k+1$ is odd.
\end{proof}
\begin{exm}[cf. {\cite[Example before Proposition 4.4]{lubotzky1999lattices}}]
    Let $k=2$, $K_0 = \mathbb{Q}_2$, and $K = \mathbb{Q}_2(\sqrt{2})$, so $r_2 = e(K/\mathbb{Q}_2) = 2$. By Theorem \ref{thm:width_formula}, $r_2/2 = 1\leq d_2\leq 1 = r_2 - 1$, so the fixed subtree in $T_K = \pgl(2,K)/\pgl(2,\mathcal{O}_K)$ of an element of order $2^k = 4$ in $\slg(2, K)$ is the $1$-neighborhood of an edge. In fact, the Eisenstein polynomial is $f(x) = x^2 - 2$, meaning $a_0 = 1$ and $a_1 = 0$. The equation for $j=1$,
    \[
    a_1 = \sqrt{a_0}\in \mathbb{F}_2
    \]
    does not hold, so our width formula implies $d_2 = \frac{r_2+j-1}{2} = 1$, in agreement with the bounds above.

    By contrast, for elements of order $8$ in $\slg(2,K)$ one has $k=3$, $\mathbb{Q}_2(z_8) = \mathbb{Q}_2(\sqrt{2}) = K$, and $r_3 = e(K/K) = 1$, so $d_3 = \frac{r_3-1}{2} = 0$ and the fixed subtree is a single edge.
\end{exm}

\section{\texorpdfstring{$N_{\gl(2,K)}(C_{2^k}.2)$}{NGL(2,K)(C2k.2)}-action on the Fixed Tree Boundary}\label{append:b}
\begin{prop}\label{prop:normalizer_not_transitive}
    In the setting of Proposition \ref{prop:dicyclic_fixed}, each $N_{\gl(2,K)}(A)$-orbit of $\partial T_A$ has at most $2$ vertices if $k\geq 3$. If $k=2$, then each $N_{\gl(2,K)}(A)$-orbit of $\partial T_A$ has at most $6$ vertices.
\end{prop}
\begin{proof}
    Fix the dicyclic group presentation $A = \langle g,\sigma\mid g^{2^{k-1}} = \sigma^2 = -I,\ \sigma g\sigma^{-1} = g^{-1}\rangle$. 
    
    If $k\geq 3$, then for any $t\in N_{\gl(2,K)}(A)$, since the $t$-adjoint map preserves the eigenvalues, we have
    \[
    tgt^{-1} \in \{g,g^{-1}\},\ t\sigma t^{-1} = g^j\sigma,\ \exists j\in\mathbb{Z}.
    \]
    On the other hand, the $\sigma$-adjoint map sends $(g,\sigma)$ to $(g^{-1},\sigma)$ while the $g$-adjoint map sends it to $(g,g^2\sigma)$. Therefore, the image of the map
    \[
    N_{\gl(2,K)}(A)\to \operatorname{Out}(A)
    \]
    is contained in a cyclic group of order $2$, generated by $(g,\sigma)\mapsto (g,g\sigma)$. Since its kernel is $K^\times A$ and fixes $T_A$ pointwise, each $N_{\gl(2,K)}(A)$-orbit of $\partial T_A$ contains at most two vertices.

    If $k = 2$, then $A$ is isomorphic to the quaternion group $Q_8$, allowing a larger outer automorphism group $\operatorname{Out}(Q_8)\cong S_3$, generated by the permutations of $\{g,\sigma,g\sigma\}$. For the same reason, each $N_{\gl(2,K)}(A)$-orbit of $\partial T_A$ contains at most $6$ vertices.
\end{proof}
\begin{rmk}
    If $k\geq 3$, then $t = I+g$ satisfies $tgt^{-1} = g$ and $t\sigma t^{-1} = g\sigma$. For $k=2$, the elements $[I+g]$, $[I+\sigma]\in N_{\gl(2,K)}(A)/(K^\times A)$ serve as the generators of this copy of $S_3$.
\end{rmk}
\bibliography{reference}{}
\bibliographystyle{plain}
		
\end{document}